\documentclass[mathpazo]{csam}

\usepackage{dmhm}

\begin{document}

\title{Distributed-memory \HMat{} Algebra I:\\
    Data distribution and matrix-vector multiplication}

\author[Y. Li et~al.]{
    Yingzhou Li\affil{1}\comma\corrauth,
    Jack Poulson\affil{2}~and
    Lexing Ying\affil{3}}

\address{
    \affilnum{1}\ Department of Mathematics, Duke University, Durham, NC
    27708, USA. \\
    \affilnum{2}\ Hodge Star, Toronto, Canada. \\
    \affilnum{3}\ Department of Mathematics and ICME, Stanford
    University, Stanford, CA 94305, USA. \\
}

\email{
    {\tt yingzhou.li@duke.edu} (Y.~Li),
    {\tt jack@hodgestar.com} (J.~Poulson),
    {\tt lexing@stanford.edu} (L.~Ying)
}

\begin{abstract}
We introduce a data distribution scheme for \HMats{} and a
distributed-memory algorithm for \HMat{}-vector multiplication. Our data
distribution scheme avoids an expensive $\Omega(P^2)$ scheduling procedure
used in previous work, where $P$ is the number of processes, while data
balancing is well-preserved. Based on the data distribution, our
distributed-memory algorithm evenly distributes all computations among $P$
processes and adopts a novel tree-communication algorithm to reduce the
latency cost. The overall complexity of our algorithm is $O\Big(\frac{N \log
N}{P} + \alpha \log P + \beta \log^2 P \Big)$ for \HMats{} under weak
admissibility condition, where $N$ is the matrix size, $\alpha$ denotes the
latency, and $\beta$ denotes the inverse bandwidth. Numerically, our
algorithm is applied to address both two- and three-dimensional problems of
various sizes among various numbers of processes. On thousands of processes,
good parallel efficiency is still observed.
\end{abstract}



\maketitle

\section{Introduction}

For linear elliptic partial differential equations, the blocks of both
forward and backward operators, when restricted to non-overlapping domains,
are numerically low-rank~\cite{Bebendorf2003}. Hence both operators can be
represented in a data sparse form. Many fast algorithms benefit from this
low-rank property and apply these operators in quasi-linear scaling. Such
fast algorithms include but not limit to tree-code~\cite{Barnes1986,
Salmon1994}, fast multipole method (FMM)~\cite{Anderson1992, Cheng1999,
Greengard1987, Greengard1990, Greengard1997, Rokhlin1985, Singh1993,
Ying2003}, panel clustering method~\cite{Hackbusch1999}, etc. The low-rank
structures in these fast algorithms are revealed via various interpolation
techniques such as: pole expansion, Chebyshev interpolation, equivalent
interaction, etc~\cite{Rokhlin1985, Greengard1987, Ying2003, Fong2009}.

In contrast to approximating the application of operators, another group of
research focuses on approximating operators directly in compressed matrix
forms. As one of the earliest members in this group,
\HMat{}~\cite{Bebendorf2007, Bebendorf2008, Bebendorf2003, Grasedyck2003,
Hackbusch1999, Hackbusch2000a, Izadi2012, Izadi2012z, Kriemann2005}
hierarchically compresses operators restricted to far-range interactions by
low-rank matrices. The memory cost and matrix-vector multiplication
complexity are quasi-linear with respect to the degrees of freedom (DOFs) in
the problem. Shortly after introducing \HMat{}, \citet{Hackbusch2000b} again
introduced $\mathcal{H}^2$-matrix, which uses nested low-rank bases to
further reduce the memory cost and multiplication complexity down to linear.
Related to the fast algorithms above, \HMat{} and $\mathcal{H}^2$-matrix can
be viewed as algebraic versions of tree code and FMM respectively. But they
are more flexible in choosing different admissibility conditions and
low-rank compression techniques, which are related to general advantages of
algebraic representations.

Developments in the \HMat{} group and extensions beyond the group are
explored in the past decade. Hierarchical off-diagonal low-rank matrix
(HOLDER)~\cite{Aminfar2016} and hierarchical semi-separable matrix
(HSS)~\cite{Xia2010a} are two popular hierarchical matrices with the
simplest admissibility condition, i.e., weak admissibility condition.
Different from hierarchical matrices, recursive skeletonization
factorization (RS)~\cite{Minden2017} and hierarchical interpolative
factorization (HIF)~\cite{Ho2016, Ho2016a} introduce separators in the
domain partition and compress the operator as products of sparse matrices.
The partition and factorization in RS and HIF are in the similar spirit as
that in multifrontal method~\cite{Duff1986, Amestoy2011} and superLU
method~\cite{Li2011}, while extra low-rank approximations are introduced to
compress the interactions within frontals. Other algebraic representations
include block low-rank approximation~\cite{Xing2018}, block basis
factorization~\cite{Wang2019a}, etc. The benefits of algebraic
representations over analytical fast algorithms come in two folds: 1)
numerical low-rank approximation is more effective than interpolation; 2)
matrix factorization and inversion become feasible. We emphases that these
algebraic representations are not only valid for linear elliptic operators,
but also valid for operators associated with low-to-medium frequency
Helmholtz equations and radial basis function kernel matrices. When
operators admit high-frequency property, the low-rank structure appears in a
very different way comparing to that in all aforementioned fast algorithms,
and are also well-studied by the community~\cite{Engquist2007, Candes2009,
Engquist2009, ONeil2010, Benson2014, Li2015a, Li2015b, Li2017, Li2018}.

Many of these fast algorithms and algebraic representations have been
parallelized on either shared-memory or distributed-memory setting to be
applicable to practical problems of interest~\cite{Greengard1990,
Salmon1994, Singh1993, Warren1993, Ying2003, Wang2013, Ghysels2016,
Rouet2016, Li2017a, Chen2018, Wang2019b, Takahashi2020}. Here we focus on
the parallelization of \HMat{}. \citet{Kriemann2005,Kriemann2013}
implemented a shared-memory parallel \HMat{} using a block-wise
distribution, i.e., each block is assigned to a single process. Processes
assigned to blocks near root level are responsible for computations of
complexity linear in $N$, where $N$ is the total DOFs. Hence the speedup of
such a parallelization scheme is theoretically upper bounded by $O(\log N)$
and limited in practice up to $16$ processes. \citet{Izadi2012z, Izadi2012}
published detailed algorithms for \HMat{} addition, matrix-vector
multiplication, matrix-matrix multiplication and matrix inversion under
distributed-memory setting. In \cite{Izadi2012z, Izadi2012}, the data of
\HMat{} are evenly distributed among all processes according to their global
matrix indices, which is similar to our data distribution for
one-dimensional problems with uniform discretization but different from ours
for other setups. The computations in \cite{Izadi2012z, Izadi2012} are
distributed under task-based parallelization, whose the scheduling part
costs $\Omega (P^2)$ operations on $P$ processes. According to numerical
results therein, good parallel efficiency is limited up to $16$ processes.

\subsection{Contribution}

In this paper, we first propose a balanced data distribution scheme for
\HMats{} based on the underlying domain geometry\footnote{When the domain
geometry of the problem is not available and only the graph connectivity of
the problem is known, our data distribution scheme can be extended to use
the hierarchical partition of the graph instead.}. In \HMat{}, the domain is
usually hierarchically partitioned and then organized in a domain tree
structure. In order to avoid any expensive scheduling procedure, our
processes are also organized in a tree structure in correspondence to that
of the hierarchical domain partition. Each process then owns a unique piece
of the domain and also own the associated data in \HMat{}. Following such a
data distribution, all data in \HMat{} are evenly distributed among all
processes. For a \HMat{} of size $N$ distributed on $P$ processes, the
memory cost is $O\Big(\frac{N\log N}{P} \Big)$ on each process. Our data
distribution scheme is scalable up to $P = O(N)$ processes.

Building on top of our data distribution, a distributed-memory parallel
algorithm is proposed to conduct the \HMat{}-vector multiplication. Our
parallel algorithm consists of several parts: a computation part, three
consecutive communication parts, and another computation part. When the
input and output vectors are distributed according to the tree structure of
processes, both computation parts are communication-free. Then a novel data
communication scheme, known as the tree-communication, is introduced to
significantly reduce costs in two of the communication parts. The remaining
communication part consists of a constant number of point-to-point
communication on each process. Mainly due to the process organization and
the tree-communication scheme, the expensive scheduling procedure is totally
avoided throughout our algorithm. The overall computational and
communication complexities, then, are $O\Big(\frac{N \log N}{P}\Big)$ and
$O\Big(\alpha \log P + \beta \Big(\log^2 P + \log \frac{N}{P} +
\big(\frac{N}{P} \big)^{\frac{d-1}{d}} \Big) \Big)$\footnote{This is
complexity for \HMats{} under standard admissibility conditions and an upper
for \HMats{} under weak admissibility condition. } respectively, where $d$
is the dimension of the problem, $\alpha$ denotes the message latency, and
$\beta$ denotes the inverse bandwidth.

Finally, the parallel algorithm is applied to two-dimensional and
three-dimensional problems of sizes varying from a few thousands to a
quarter billion on massive number of processes. The parallel scaling is
still found to be near-ideal on computational resources available to us,
up to a few thousands processes. In all cases, our \HMat{}-vector
multiplications are completed within a few seconds.

\subsection{Organization}

The rest of the paper is organized as follows. In
Section~\ref{sec:preliminary}, we revisit \HMat{} together with
admissibility conditions. Section~\ref{sec:hmat-data-distribution}
introduces our balanced data distribution scheme. The distributed-memory
\HMat{}-vector multiplication algorithm is detailed in
Section~\ref{sec:hmat-vector-multiplication}.
Section~\ref{sec:numerical-results} presents numerical results for
two-dimensional and three-dimensional problems of various sizes. Finally, we
conclude the paper in Section~\ref{sec:conclusion} together with some
discussion on future work.

\section{Preliminary}
\label{sec:preliminary}

In this section, we first review the definition and the structure of
\HMat{}. Then the \HMat{}-vector multiplication follows in a
straightforward way.

Let us assume that $\calK(t,k)$ is a kernel satisfying the hierarchical
low-rank property as in tree code or \HMat{}. Then applying Nystr{\"o}m
discretization to the integral equation,
\begin{equation} \label{eq:intop}
    u(t) = \int_\Omega \calK(t,s) f(s) \dif s, \quad \text{for } t\in
    \Omega,
\end{equation}
results a matrix-vector multiplication, and the matrix therein can
be approximated by an \HMat{}. Throughout the rest paper, we use
the concepts of a domain and the Nystr{\"o}m discretization points
in the domain interchangably. For example, a matrix restricted to
$\Omega_1 \times \Omega_2$ means that the matrix restricted to the
row and column indices corresponding to the discretization points
in $\Omega_1$ and $\Omega_2$ respectively. In \eqref{eq:intop},
the operator maps from the domain $\Omega$ to itself. In practice,
\HMat{} can also be used to approximate operators mapping from one
domain to another and the rest of the paper can be extended to such
a setting with a minor update on domain notations. To simplify our
presentation, we limit ourselves to the self mapping case.

\begin{figure}[ht]
    \centering
    \includegraphics[width=0.8\textwidth]{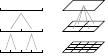}
    \caption{Hierarchical partition of ideal one-dimensional domain
    (left) and two-dimensional domain (right). Gray lines indicate
    connections between parent domains and their child subdomains.}
    \label{fig:domain-partition}
\end{figure}

In the above setting, the structure of the \HMat{} fundamentally relies on
the hierarchical partition of the domain $\Omega$, which is defined as
follows:
\begin{definition}[\bfemph{Domain tree}]
    A tree $\treeO = (\calV_\Omega, \calE_\Omega)$ with the vertex set
    $\calV_\Omega$ and the edge set $\calE_\Omega$ is called a domain tree
    of $\Omega$ if the following conditions hold:
    \begin{enumerate}
        \item All nodes in $\treeO$ are subdomains of $\Omega$;
        \item The set of children of a domain $\omega \in
        \calV_\Omega$, denoted as $\calC(\omega) = \set{\nu \in
        \calV_\Omega \mid \exists (\omega, \nu) \in \calE_\Omega}$, is
        either empty or a partition of $\omega$;
        \item $\Omega \in \calV_\Omega$ is the root of $\treeO$.
    \end{enumerate}
\end{definition}

When a (quasi-)uniform discretization of a regular $d$-dimensional domain
$\Omega = [0,1]^d$ is considered, the domain tree is constructed via
applying a $2^d$ uniform partition recursively. Such domains are later
referred as ideal $d$-dimensional domains. Figure~\ref{fig:domain-partition}
illustrates two domain tree associated with an ideal one-dimensional domain
and an ideal two-dimensional domain.

The low-rank submatrices in an \HMat{} are determined by admissibility
conditions. There are many different admissibility conditions leading to
different \HMat{} structures. Here we introduce two of them: weak
admissibility condition and standard admissibility condition.

\begin{figure}[ht]
    \centering
    \includegraphics[width=0.8\textwidth]{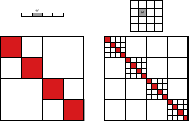}
    \caption{Weak admissibility condition and the corresponding \HMats{}
    for ideal one-dimensional (left) and two-dimensional (right) domains.
    First row shows domain partitions and gray blocks are non-admissible
    domains to $\omega$. The second row shows the corresponding \HMats{}
    with red submatrices being dense and white ones being low-rank.}
    \label{fig:weak-ad-cond}
\end{figure}

\begin{figure}[ht]
    \centering
    \includegraphics[width=0.8\textwidth]{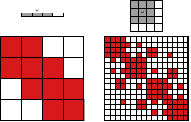}
    \caption{Standard admissibility condition and the corresponding \HMats{}
    for ideal one-dimensional (left) and two-dimensional (right) domains.}
    \label{fig:standard-ad-cond}
\end{figure}

\begin{definition}[\bfemph{Weak admissibility condition}]
    Two domains, $\omega$ and $\nu$, are weakly admissible if $\omega \cap
    \nu = \emptyset$.
\end{definition}

\begin{definition}[\bfemph{Standard admissibility condition}]
    Two domains, $\omega$ and $\nu$, are standard admissible if
    \begin{equation} \label{eq:standard-ad-cond}
        \min \big(\mathrm{diam}(\omega), \mathrm{diam}(\nu) \big) \leq
        \rho \mathrm{dist} (\omega, \nu),
    \end{equation}
    where $\mathrm{diam}(\omega)$ is the diameter of $\omega$,
    $\mathrm{dist}(\omega, \nu)$ is the distance between two domains, and
    $\rho$ is a constant adjusting the size of buffer zone.
\end{definition}

Weak admissibility condition is the simplest admissibility condition used in
practice and leads to the simplest \HMat{} structure. While standard
admissibility condition is more complicated, but widely used in many fast
algorithms~\cite{Barnes1986, Greengard1987}. Importantly, for linear
elliptic differential operators with $L_\infty$ coefficients discretized by
a local basis set, both the forward differential operator and its inverse
can be well-approximated by \HMat{} under standard admissibility condition.
Throughout this paper, we adopt $\rho \equiv \sqrt{d}$.
Figure~\ref{fig:weak-ad-cond} and Figure~\ref{fig:standard-ad-cond} shows
the weak admissibility condition and the standard admissibility condition
respectively for both ideal one-dimensional and two-dimensional domains. One
more popular admissibility condition, known as strong admissibility
condition~\cite{Bebendorf2008, Minden2017}, simply replaces the ``$\min$''
in \eqref{eq:standard-ad-cond} by ``$\max$''. Strong admissibility condition
and standard admissibility condition are the same on ideal domains.

Based on the domain tree and admissibility conditions, we are ready to
precisely define \HMat{} as an approximation of $\calK$ mapping from
$\Omega_s = \Omega$ to $\Omega_t = \Omega$, i.e., an approximation of
$\calK_{\Omega_t \times \Omega_s} = \calK\rvert_{\Omega_t \times \Omega_s}$.
Domain $\Omega_t$ and $\Omega_s$ are called \bfemph{the target domain} and
\bfemph{the source domain} respectively, where the target domain is
associated with row indices and the source domain is associated with column
indices.

\begin{definition}[\HMat{}] \label{def:hmat}
    Assume that $\calK$ maps vectors defined on source domain $\Omega_s$
    to vectors defined on target domain $\Omega_t$. $\calK$ is an \HMat{}
    with rank $r$ and domain trees $\treeOt$ and $\treeOs$ if the
    following conditions hold in order: for each child subdomain pairs of
    $\Omega_t \times \Omega_s$, i.e., $\omega_t \times \omega_s \in
    \calC(\Omega_t) \times \calC(\Omega_s)$,
    \begin{enumerate}
        \item if $\calC(\omega_t) = \emptyset$ or $\calC(\omega_s) =
        \emptyset$, then $\calK_{\omega_t \times \omega_s}$ is \bfemph{a
        dense matrix} $D_{\omega_t \times \omega_s}$; else
        \item if $\omega_t$ and $\omega_s$ are admissible, then
        $\calK_{\omega_t \times \omega_s}$ is \bfemph{a low-rank matrix}
        with rank $r$, i.e., $\calK_{\omega_t \times \omega_s} =
        U_{\omega_t \times \omega_s} V_{\omega_t \times \omega_s}^\top$
        for $U_{\omega_t \times \omega_s} \in \bbR^{\lvert \omega_t \rvert
        \times r}$, $V_{\omega_t \times \omega_s} \in \bbR^{\lvert
        \omega_s \rvert \times r}$, and $\lvert \cdot \rvert$ denotes the
        DOFs in the domain; otherwise
        \item $\calK_{\omega_t \times \omega_s}$ is \bfemph{an \HMat{}}
        with rank $r$ and domain trees $\bbT_{\omega_t}$ and
        $\bbT_{\omega_s}$.
    \end{enumerate}
\end{definition}

In the above definition, three conditions must be checked in the given
order. The third condition defines the hierarchical structure of \HMat{}.

In order to further clarify the definition, we walk readers through the
ideal two-dimensional domain case under weak admissibility condition. We
begin with $\calK$ mapping from $\Omega_s = \Omega = [0,1]^2$ to $\Omega_t =
\Omega = [0,1]^2$. There are 16 child subdomain pairs in $\calC(\Omega_t)
\times \calC(\Omega_s)$. Among all 16 pairs, all child domains have their
child domains. Hence the first condition in Definition~\ref{def:hmat} fails
for all pairs. We then check the weak admissibility condition. There are 12
out of 16 pairs are admissible, i.e., all non-overlapping domain pairs on
the first level as in Figure~\ref{fig:domain-partition}. Therefore, there
are 12 off-diagonal submatrices are low-rank, which are denoted by the big
white blocks in Figure~\ref{fig:weak-ad-cond} (right). For the rest 4 child
subdomain pairs, they are \HMats{} of them own. We can continue this process
until the leaf level of the domain tree and resolve the entire \HMat{} as in
Figure~\ref{fig:weak-ad-cond} (right). The \HMat{} under standard
admissibility condition is much more complicated.
Figure~\ref{fig:standard-ad-cond} depicts the \HMats{} with the same domain
and domain tree as that in Figure~\ref{fig:weak-ad-cond} but under the
standard admissibility condition instead.

The \HMat{}-vector multiplication can be processed efficiently as long as we
can read from the input vector and write to the output vector restricting to
subdomains. We denote the \HMat{} as $\calK$ mapping from $\Omega$ to
$\Omega$ and the \HMat{}-vector multiplication as,
\begin{equation*}
    y = \calK x,
\end{equation*}
where both $x$ and $y$ are vectors defined on $\Omega$. We first
initialize the output vector $y$ as a zero vector. Then, we traverse all
submatrices in $\calK$ that contains data, i.e., dense submatrices and
low-rank submatrices. For any such submatrix, denoted as $\calK_{\omega_t
\times \omega_s}$, we conduct the matrix-vector multiplication and add the
results to the output vector,
\begin{equation} \label{eq:hmat-vec-add}
    y_{\omega_t} = y_{\omega_t} + \calK_{\omega_t
    \times \omega_s} x_{\omega_s} =
    \begin{cases}
        y_{\omega_t} + D_{\omega_t \times \omega_s} x_{\omega_s} \\
        y_{\omega_t} + U_{\omega_t \times \omega_s} (V_{\omega_t \times
        \omega_s}^\top x_{\omega_s}) \\
    \end{cases},
\end{equation}
where $y_{\omega_t}$ and $x_{\omega_s}$ denote the vector restricted to
domain $\omega_t$ and $\omega_s$ respectively. When
\eqref{eq:hmat-vec-add} is completed for all submatrices in $\calK$, the
vector $y$ is already the final \HMat{}-vector multiplication result.

As shown in many previous work~\cite{Bebendorf2008,Hackbusch1999}, for the
regular domain with almost uniformly distributed discretization points, the
memory cost for the \HMat{} is $O(r N \log N)$, where $N$ is the total DOFs
and $r$ is the numerical rank. The \HMat{}-vector multiplication can be
achieved in $O(r N \log N)$ operations. Here we mainly reviewed the
structure and matrix-vector multiplication of \HMat{}. The construction
algorithms of \HMat{}~\cite{Bebendorf2008, Lin2011} as well as other
algebraic operations such as: matrix-matrix multiplication, matrix
factorization, etc., have been extensively studied in the literature, which
are beyond the scope of this paper and we omit the detailed discussion.

\section{\HMat{} Data Distribution}
\label{sec:hmat-data-distribution}

The organization of processes and the associated data distribution avoid
expensive parallel scheduling procedure as in~\cite{Izadi2012,
Izadi2012z}. In Section~\ref{sec:process-organization}, we first explain
our hierarchical organization of processes. Then in
Section~\ref{sec:data-distribution}, the data distribution together with
the load balancing are discussed.

\subsection{Hierarchical Process Organization}
\label{sec:process-organization}

Processes are organized in correspondence with the domain tree $\treeO$. The
main idea is to assign subdomains to processes as balanced as possible
while preserving the hierarchical structure.

Let the $P$ processes be indexed from 0 to $P-1$, and $P$ be upper bounded by
the number of leaf nodes in $\treeO$~\footnote{Having more processes than
the number of leaf nodes ($O(N)$) is feasible if the later algorithm
description is slightly modified. While, such a setup is not of practical
usage. Hence we omit the detail.}. The set of all processes, denoted as
$\calP = \{0, 1, \dots, P-1\}$, is called \bfemph{the process group}. We
then traverse the domain tree to assign the process group, subgroups, or
individual processes to nodes in $\treeO$. Regarding the root node in
$\treeO$, i.e., domain $\Omega$, we assign the entire process group $\calP$
to it. From now on, we consider a general domain $\Omega^\ell$ in $\treeO$
at level $\ell$ with a general process group $\calP^\ell$ assigned. The
assignment of subgroups of $\calP^\ell$ to child subdomains of $\Omega^\ell$
obeys the following conditions:
\begin{enumerate}[label=\arabic*)]
    \item If the number of child subdomains of $\Omega^\ell$ is smaller than
    or equal to the number of processes in $\calP^\ell$, i.e., $\lvert
    \calC(\Omega^\ell) \rvert \leq \lvert \calP^\ell \rvert$, then
    $\calP^\ell$ is partitioned into $\lvert \calC(\Omega^\ell) \rvert$
    subgroups such that the number of processes in each subgroup is
    proportional to the DOFs in the corresponding child subdomain. Each
    subgroup is then assigned to the corresponding child subdomain.
    \item If the number of child subdomains of $\Omega^\ell$ is bigger than
    the number of processes in $\calP^\ell$, i.e., $\lvert
    \calC(\Omega^\ell) \rvert > \lvert \calP^\ell \rvert$, then
    $\calC(\Omega^\ell)$ are organized into $\lvert \calP^\ell \rvert$ parts
    such that the total DOFs in each part are balanced. Each process is then
    assigned to subdomains in one part.
\end{enumerate}

\begin{figure}[htp]
    \centering
    \includegraphics[width=0.8\textwidth]{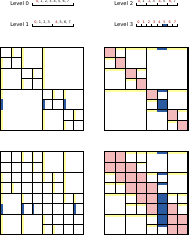}
    \caption{\HMat{} data distribution of an ideal one-dimensional domain
    with weak and standard admissibility condition on 8 processes. Top part
    is the the four-level domain tree together with its process assignment.
    Red processes are group leaders. Middle parts and bottom parts are the
    distributed \HMat{} with weak and standard admissibility condition
    respectively. Left columns are \HMat{} owned by target process groups
    and right columns are owned by source process groups. Blue blocks
    indicate data owned by process 5 whereas light red and yellow blocks
    indicate data owned by other processes.}
    \label{fig:process-tree-1d-standard}
\end{figure}

\begin{figure}[ht]
    \centering
    \includegraphics[width=0.8\textwidth]{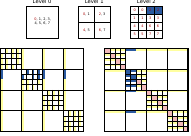}
    \caption{\HMat{} data distribution of an ideal two-dimensional domain
    with weak admissibility condition on 8 processes. Top part is the the
    three-level domain tree together with its process assignment. Red
    processes are group leaders. Bottom parts are data in \HMat{} owned by
    target process groups (left) and source process groups (right). Blue
    blocks indicate data owned by process 2 whereas light red and yellow
    blocks indicate data owned by other processes.}
    \label{fig:process-tree-2d-weak}
\end{figure}

According to the process organization strategy, a process participates and
only participates one process group at each level. When a single process is
assigned to a subdomain $\omega$ in $\treeO$, it is assigned to all
descendants of $\omega$ in $\treeO$. We can then combine all subdomains that
are singly owned by a process $p$ and denote the union as $\omega_p$. All
such unions $\{\omega_p\}_{p=0}^{P-1}$ form a balanced partition of
$\Omega$, where the balancing factor is upper bounded by twice the balancing
factor of $\treeO$. The balancing factor here is referring to the ratio of
the heaviest workload and the lightest workload among all processes.
Further, in each process group or subgroup, the process with smallest index
is called \bfemph{the group leader}, e.g., 0 is the group leader of $\calP$.
When a process is the group leader at a level $\ell$, then it is the group
leader in all descendant groups it participates. For example, process 0 is
group leaders of all process groups it participates, whose workload is the
heaviest among all processes. 

Figure~\ref{fig:process-tree-1d-standard} top show a four level domain tree
together with its process assignment of $\calP = \{0, 1, \dots, 7\}$ for an
ideal one-dimensional domain. Each process is assigned to a unique subdomain
at level 3. In addition, Figure~\ref{fig:process-tree-2d-weak} top show  a
three level domain tree together with its process assignment of $\calP =
\{0, 1, \dots, 7\}$ for an ideal two-dimensional domain. We find that
process owns two subdomains at level 2 and eight processes form a perfect
partition of the domain.

\subsection{Data Distribution and Load Balancing}
\label{sec:data-distribution}

Definition~\ref{def:hmat} explicitly shows that all data (matrix entries)
are in two types of submatrices, either dense submatrices or low-rank
submatrices. The hierarchical submatrices defined by the third conditions in
Definition~\ref{def:hmat} exist virtually for recursion purpose. Hence, we
just need to distribute dense submatrices and low-rank submatrices among
processes.

\bfemph{Low-rank submatrix.} Consider a low-rank submatrix associated with
domain pair $\Omega_t \times \Omega_s$, where $\Omega_t$ and $\Omega_s$ are
the target and source domains respectively. According to the process
organization defined in Section~\ref{sec:process-organization}, there are
two process groups assigned to $\Omega_t$ and $\Omega_s$, denoted as
$\calP_t$ and $\calP_s$ respectively. We distribute two factors in the
low-rank submatrices, $U_{\Omega_t \times \Omega_s}$ and $V_{\Omega_t \times
\Omega_s}$, to two process groups, i.e., $U_{\Omega_t \times \Omega_s}$ is
stored among $\calP_t$ and $V_{\Omega_t \times \Omega_s}$ is stored among
$\calP_s$. If either $\calP_t$ or $\calP_s$ has only one process, then the
process owns the entire matrix. Now, assume there are more than one process
in $\calP_t$. Since $U_{\Omega_t \times \Omega_s}$ is a tall and skinny
matrix, it is distributed in a block row fashion. For each process $p \in
\calP_t$, the rows corresponding to $\omega_p$ is owned by process $p$,
where $\omega_p$ is the singly owned subdomain of $p$. If there are more
than one process in $\calP_s$, then $V_{\Omega_t \times \Omega_s}$ is
distributed in the same way among processes in $\calP_s$.

\bfemph{Dense submatrix.} Consider a dense submatrix associated with domain
pair $\Omega_t \times \Omega_s$ and the corresponding process group pair
$\calP_t \times \calP_s$. There are three scenarios of the sizes of process
groups: (i) $\lvert \calP_t \rvert = \lvert \calP_s \rvert = 1$; (ii)
$\lvert \calP_t \rvert = 1$ and $\lvert \calP_s \rvert > 1$; (iii) $\lvert
\calP_t \rvert > 1$ and $\lvert \calP_s \rvert = 1$. In the first scenario,
the dense matrix $D_{\Omega_t \times \Omega_s}$ is owned by $\calP_s$. In
the second scenario, the transpose of dense matrix, $D_{\Omega_t \times
\Omega_s}^\top$, is distributed among $\calP_s$ in the same way as the
distribution of $V_{\Omega_t \times \Omega_s}$ above. In the last scenario,
the dense matrix $D_{\Omega_t \times \Omega_s}$ is distributed among
$\calP_t$ in the same way as the distribution of $U_{\Omega_t \times
\Omega_s}$ above.

\vspace{1em}
Once the data distribution strategies are applied to all submatrices, the
\HMat{} is then fully distributed among $\calP$. To further facilitate the
understanding of the overall data distribution,
Figure~\ref{fig:process-tree-1d-standard} and
Figure~\ref{fig:process-tree-2d-weak} show the data in \HMats{} for an ideal
one-dimensional domain with weak and standard admissibility condition and an
ideal two-dimensional domain with weak admissibility condition respectively.
Both \HMats{} are distributed among process group $\calP$ of size eight. In
Figure~\ref{fig:process-tree-1d-standard} and
Figure~\ref{fig:process-tree-2d-weak}, blue blocks highlight the data owned
by process 5 and process 2 respectively.

\begin{remark}
    The data distribution strategies we introduced here are suitable
    and efficient for a sequence of parallel-friendly \HMat{}
    algebraic operations, e.g., matrix-vector multiplication,
    matrix-matrix multiplication, matrix compression, matrix addition,
    etc. While some other \HMat{} algebraic operations, like \HMat{}-LU
    factorization and \HMat{}-inversion, are not parallel-friendly
    since the operations therein depends sequentially on each
    other. Our data distribution strategies work for these operations
    as well, while the efficiency is left to be further explored.
\end{remark}

\begin{remark}
    These data distribution strategies can also be easily extended to
    $\mathcal{H}^2$-matrix. The nested basis in $\mathcal{H}^2$-matrix
    can be distributed among all processes in the similar way as we
    distribute low-rank factors. While the tiny middle matrix in each
    low-rank block in $\mathcal{H}^2$-matrix could be singly owned by
    either its source or target group leader. Given such distribution
    strategies for $\mathcal{H}^2$-matrix, all its algebraic operations
    can be parallelized in an analog way as that for \HMat{}.
\end{remark}

We now discuss the load balancing of the distributed \HMat{}. As shown in
Figure~\ref{fig:process-tree-1d-standard} and
Figure~\ref{fig:process-tree-2d-weak}, the load balancing is different for
different admissibility conditions. Figure~\ref{fig:process-tree-2d-weak}
under weak admissibility condition shows an ideal load balancing whereas
Figure~\ref{fig:process-tree-1d-standard} under standard admissibility
condition shows slightly unbalanced data distribution. In the following, we
assume the domain is an ideal $d$-dimensional domain, $[0,1]^d$ with $n$
uniform discretization points on each dimension and $N = n^d$ discretization
points in total. In such an ideal case, each process own the same size of
subdomain.

Assume that the weak admissibility condition is applied. At each level on $\treeO$, any domain
$\Omega^\ell$ has the same number of admissible domains. Each process participate one domain on the
target side and another on the source side. Hence all processes own exactly the same amount of data
in low-rank submatrices at each level. For all low-rank submatrices throughout levels, data are
evenly distributed among all processes. Regarding the dense submatrices, they are all of the same
size and owned by their source processes. Since all processes own the same size domains on the
source side, and these domains have the same amount of dense submatrices, all processes own the same
amount of dense submatrix data. Overall, the data of dense submatrices and low-rank submatrices are
evenly distributed among all processes and the load balancing in this case is ideal.

While, when the standard admissibility condition is applied, the load
balancing depends on the boundary condition of the problem. If the periodic
boundary condition is adopted, the load balancing is still ideal. While, if
a non-periodic boundary condition is adopted, the data loads are different
for processes owning domains near the center and processes owning domains
near corners. Since all low-rank submatrices are evenly owned by processes
in its process groups, the load balancing factor is simply the ratio of the
numbers of low-rank submatrices for different processes, i.e., the numbers
of admissible domains. Consider level $\ell$, which is neither the first two
levels nor the last one. A center subdomain $\Omega^\ell_\text{center}$'s
parent domain has $3^d$ non-admissible neighbor domains, each of which is
partitioned into $2^d$ subdomains at level $\ell$. Excluding non-admissible
subdomains of $\Omega^\ell_\text{center}$, there are $3^d \cdot 2^d - 3^d$
admissible subdomains of $\Omega^\ell_\text{center}$. However, a corner
subdomain $\Omega^\ell_\text{corner}$'s parent domain is also a corner
domain and has $2^d$ non-admissible neighbor domains. Through the similar
calculation, $\Omega^\ell_\text{corner}$ has $2^d \cdot 2^d - 2^d$
admissible subdomains at level $\ell$. Hence the load balancing factor is
$\frac{3^d}{2^d}$. Such a factor also holds to the load balancing of dense
submatrices. Overall, asymptotically as $N$ goes to infinity, the load
balancing factor for distributed \HMat{} under standard admissibility
condition and non-periodic boundary condition is upper bounded by
$\big(\frac{3}{2} \big)^d$. Since this factor is independent of both $N$ and
$P$, we still regard our data distribution in this case as a balanced one.

\section{Distributed-memory \HMat{}-vector Multiplication}
\label{sec:hmat-vector-multiplication}

\HMat{}-vector multiplication is the fundamental operation in \HMat{}
algebra and reveals the value of \HMat{} as a fast algorithm. Further, it is
also one of basic operations involved in other \HMat{} algebraic operations,
including, matrix-matrix multiplication, matrix compression, matrix
factorization, and matrix inversion. As briefly reviewed in
Section~\ref{sec:preliminary}, the sequential \HMat{}-vector multiplication
is as simple as looping over all low-rank and dense submatrices, multiplying
the submatrix to the input vector restricted to the source domain, and
adding the result to the output vector restricted to the target domain.
However, the distributed-memory version is much more complicated. Based on
the data distribution as in Section~\ref{sec:hmat-data-distribution}, we
present the distributed-memory \HMat{}-vector multiplication algorithm in
this section followed by its complexity analysis.

\subsection{Algorithm}
\label{sec:algorithm}

Distributed-memory \HMat{}-vector multiplication algorithm mainly consists
of the following five steps:
\begin{enumerate}[label=Step \arabic*., leftmargin=3.5\parindent]
    \item Source side local computation;
    \item Tree-reduction on source process tree;
    \item Data transfer from source to target;
    \item Tree-broadcast on target process tree;
    \item Target side local computation.
\end{enumerate}
Among these five steps, Steps 1 and 5 only involve computations and are
communication-free whereas Steps 2, 3, and 4 focus on efficient
communication under our data distribution and process organization. We will
elaborate five steps in detail one-by-one. Throughout the following
description, we assume the input vector $x$ is already distributed in the
block row fashion among process group $\calP$. More precisely, for any
process $p \in \calP$, it owns $x_{\omega_p} = x \rvert_{\omega_p}$ for
$\omega_p$ being $p$'s singly owned domain. The output vector $y$ will be
distributed exactly in the same way as $x$.

\subsubsection{Source Side Local Computation}
\label{sec:step1}

The source side local computation goes through all submatrices containing
data, i.e., low-rank submatrices and dense submatrices, and conducts all
communication-free calculations. We now describe specific operations for
submatrices of different types.

\bfemph{Low-rank submatrix.} Consider a low-rank submatrix associated with
$\Omega_t \times \Omega_s$ and process groups $\calP_t \times \calP_s$. The
explicit block form of $V_{\Omega_t \times \Omega_s}$ and $x_{\Omega_s}$
admit,
\begin{equation}
    V_{\Omega_t \times \Omega_s} =
    \begin{pmatrix}
        v_{p_0} \\
        \vdots \\
        v_{p_{\lvert \calP_s \rvert -1}}
    \end{pmatrix},
    \quad
    x_{\Omega_s} =
    \begin{pmatrix}
        x_{p_0} \\
        \vdots \\
        x_{p_{\lvert \calP_s \rvert -1}}
    \end{pmatrix},
\end{equation}
where $p_i \in \calP_s$, $v_{p_i}$ and $x_{p_i}$ are stored on process
$p_i$. We aim to compute the product of $V_{\Omega_t \times \Omega_s}$ and
$x_{\Omega_s}$ as,
\begin{equation} \label{eq:prod-vx}
    V_{\Omega_t \times \Omega_s}^\top x_{\Omega_s} = \sum_{i=0}^{\lvert
    \calP_s \rvert - 1} v_{p_i}^\top x_{p_i},
\end{equation}
where the summation over $i$ requires communication since $v_{p_i}^\top
x_{p_i}$ are owned by different processes for different $i$. Hence, in this
step, we only compute
\begin{equation}
    z_\text{local} = v_{p_i}^\top x_{p_i}
\end{equation}
on process $p_i$ without conducting any communication. The communication for
the summation over $i$ in \eqref{eq:prod-vx} is postponed until the next
step.

\bfemph{Dense submatrix.} Consider a dense submatrix associated with
$\Omega_t \times \Omega_s$ and $\calP_t \times \calP_s$. When there are more
than one process in the target process group, i.e., $\lvert \calP_t \rvert >
1$, the data in this submatrix are owned by the target process group. No
local computation is needed and we assign $z_\text{local} = x_{\Omega_s}$
for later communications. When there is only one process in the target
process group, i.e., $\lvert \calP_t \rvert = 1$, the data are distributed
among the source process group as,
\begin{equation}
    D_{\Omega_t \times \Omega_s} =
    \begin{pmatrix}
        d_{p_0} & \cdots & d_{p_{\lvert \calP_s \rvert - 1}}
    \end{pmatrix},
    \quad
    x_{\Omega_s} =
    \begin{pmatrix}
        x_{p_0} \\
        \vdots \\
        x_{p_{\lvert \calP_s \rvert -1}}
    \end{pmatrix},
\end{equation}
for $p_i \in \calP_s$ and $\lvert \calP_s \rvert \geq 1$. Similar to the
low-rank submatrix case, we aim to compute
\begin{equation} \label{eq:prod-dx}
    D_{\Omega_t \times \Omega_s} x_{\Omega_s} = \sum_{i=0}^{\lvert
    \calP_s \rvert - 1} d_{p_i} x_{p_i}.
\end{equation}
Instead, we only conduct local computation in this step, $z_\text{local} =
d_{p_i} x_{p_i}$, on each process $p_i \in \calP_s$ without communication.

\subsubsection{Tree-reduction on Source Process Tree}
\label{sec:step2}

This step implements the communication required summations in \eqref{eq:prod-vx} and
\eqref{eq:prod-dx}. Na{\"i}vely, we can perform many MPI reductions~\footnote{We refer to
  ``MPI\textunderscore Reduce'' with addition operation as the reduction throughout this paper.},
one for each submatrices and reduce the summation results to their group leaders.  However, such a
na{\"i}ve reduction strategy requires many more messages than the {\em tree-reduction} to be
introduced below, which benefits most from the hierarchical organization of both the \HMat{} and
processes.

The preliminary step in tree-reduction is to collect and pack local results
that require communication in \eqref{eq:prod-vx} and \eqref{eq:prod-dx}. For
each process, we visit the \HMat{} level by level from root to leaf. At each
level, each process participates and only participates in one process group.
Hence, local results are about to be reduced to the same group leader and
are packed together in an array in the same ordering. Across levels, we
concatenate packed local results together until one level before the level
where process group has only one process. We denote the maximum number of
such levels as $L_P$.

Then a sequence of reductions are conducted from level $L_P$ backward to the
root level. At level $L_P$, all processes reduce the entire concatenated
array to their own group leaders at this level. Group leaders at level $L_P$
then have already collected their group members' contributions to summations
from root level to level $L_P-1$. Hence, those non-leader group members at
level $L_P$ no longer participate the rest communications in this step. At a
following level $\ell = L_P - 1, L_P - 2, \dots, 1$, the participating
processes are those group leaders at level $\ell+1$. They reduce their
concatenated array from level 1 to level $\ell$ (with contributions from
their own group members) to their own group leaders at level $\ell$. When
all reductions are completed, all group leaders own the summations
\eqref{eq:prod-vx} and \eqref{eq:prod-dx} of their groups. Slightly abuse of
notation, we still denote these summation results as $z_\text{local}$.

\begin{figure}[htp]
    \centering
    \includegraphics[width=0.8\textwidth]{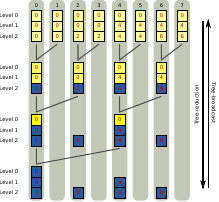}
    \caption{Tree-communication flowchart. Tree-reduction and tree-broadcast
    flow from top to bottom and from bottom to top respectively. Different
    columns with gray background are the concatenated arrays owned by
    different processes. Each cubic is the packed data on the corresponding
    level and the number in the cubic indicates its group leader. For
    tree-reduction, yellow cubics are local data to be reduced to their
    group leaders and summed together whereas blue cubics are the final
    summation results owned only by group leaders. As shown in the figure,
    only yellow cubics and their owner processes participate the reduction
    communications. For tree-broadcast, blue cubics are original packed data
    to be sent to group members. Yellow cubics are packed data been
    broadcasted. Light yellow cubics are final broadcasted data.}
    \label{fig:tree-communication-1d-weak}
\end{figure}

\begin{remark} \label{rmk:skipping}
    When the domain and discretization are far from balanced ones, the
    process tree is also not balanced. Hence, it is possible that at some
    level $\ell < L_P$, a process is the process group of its own. In this
    case, such a process do not need to participate the reduction at level
    $\ell$ or lower. We do not exclude such cases from our description
    above, but do exclude them from our implementation.
\end{remark}

Figure~\ref{fig:tree-communication-1d-weak} depicts the flow of a
tree-reduction for an ideal one-dimensional domain distributed evenly on 8
processes. Although there is no communication-required data on the root
level in \HMat{}-vector multiplication, we still include data cubics on
level 0 in the figure to demonstrate the idea and show the extendability of
the tree-reduction to more than two levels.

\subsubsection{Data Transfer from Source to Target}
\label{sec:step3}

After the previous step, all local data, $z_\text{local}$, are stored on
their own group leaders on the source side. In order to finish the
computation, local data should be sent to the processes in the target group.
To better benefit from the hierarchical structure, we accomplish the
communication in this and next steps. In this step, local data will be sent
from the source group leaders to the corresponding target group leaders.
Then the next step is responsible for broadcasting local data to the
processes in target groups.

Given a pair of target and source group leaders, $p_t$ and $p_s$, they
could be the group leaders of many submatrices. Hence process $p_s$ first
packs local data in all those submatrices and then send them in one
message to process $p_t$. After process $p_t$ received the packed local
data, it then unpacks the data to submatrices.

\begin{remark} \label{rmk:step3}
    We emphasize that a process only participates at most $O(\log P)$
    number of group leader pairs. Let us consider process 0 as the source
    group leader, which acts most frequently as the source group leader
    among all processes. As we mentioned before, each process only
    participates one process group on each level of the process tree.
    Process 0 is then the group leaders of one process group on each
    level, which adds to $O(\log P)$ groups. A source process group on
    each level only interacts with a constant number of target process
    groups, where the constant depends on the admissibility condition.
    Hence process 0 is paired with a constant number of target group
    leaders at each level. Summing all levels together, process 0 is
    paired with $O(\log P)$ target group leaders.
\end{remark}

\subsubsection{Tree-broadcast on Target Process Tree}
\label{sec:step4}

Consider a low-rank submatrices associated with $\Omega_t \times \Omega_s$
with process groups $\calP_t \times \calP_s$ as an example. The matrix
vector multiplication admits,
\begin{equation} \label{eq:prod-uvx}
    U_{\Omega_t \times \Omega_s} V_{\Omega_t \times \Omega_s}^\top
    x_{\Omega_s} =
    \begin{pmatrix}
        u_{p_0} \big( V_{\Omega_t \times \Omega_s}^\top x_{\Omega_s}
        \big) \\
        \vdots \\
        u_{p_{\lvert \calP_t \rvert - 1}}
        \big( V_{\Omega_t \times \Omega_s}^\top x_{\Omega_s}
        \big) \\
    \end{pmatrix}
    =
    \begin{pmatrix}
        u_{p_0} z_\text{local} \\
        \vdots \\
        u_{p_{\lvert \calP_t \rvert - 1}} z_\text{local}
    \end{pmatrix},
\end{equation}
where $p_i \in \calP_t$ and $z_\text{local}$ is the summation in
\eqref{eq:prod-vx}. After the previous step, in each submatrices,
$z_\text{local}$ is owned by the target group leaders. Hence, in order to
conduct the product of $u_{p_i} z_\text{local}$ as in \eqref{eq:prod-uvx},
$z_\text{local}$ needs to be shared with all target group members. A similar
equation can be written down for dense submatrices with target process
groups of size greater than one. In this step, we hierarchically broadcast
the local data $z_\text{local}$ from the group leaders to the group members
together and name it as {\em tree-broadcast}, which is the reverse procedure
of tree-reduction.

Similar to tree-reduction, we first collect and pack local results that
require communication. For each group leader, we visit the \HMat{} level by
level from root to leaf. At each level, local results that are about to be
broadcasted to the same group are packed together in an array. Across
levels, we concatenate packed local results together until level $L_P$.

Then a sequence of broadcasts are executed from the first level forward to
level $L_P$. At a level $\ell = 1, \dots, L_P - 1$, the group leaders
broadcast their array from level 1 to level $\ell$ to those subgroup leaders
at level $\ell+1$. Subgroup leaders then concatenate the received array
together with their own packed array. Once the concatenating procedure is
accomplished, we move on to the next level. Finally, at level $L_P$, group
leaders broadcast their entire array to all their group members. All
processes in target process group, in the end, received all needed local
data for each submatrices they participated. 

Similar level skipping for the unbalanced target process tree can be done
for tree-broadcast as that for tree-reduction in Remark~\ref{rmk:skipping}.
Figure~\ref{fig:tree-communication-1d-weak} illustrates a tree-broadcast
procedure for an ideal one-dimensional domain distributed on 8 processes.

\subsubsection{Target Side Local Computation}
\label{sec:step5}

The target side local computation goes through all low-rank and dense
submatrices and conducts aggregation of the product results onto output
vector $y$. Here we assume the output vector $y$ is initialized to be all
zero. We describe operations for different types of submatrices.

\bfemph{Low-rank submatrix.} Consider a low-rank submatrix associated with
$\Omega_t \times \Omega_s$ and process groups $\calP_t \times \calP_s$. As
shown in \eqref{eq:prod-uvx}, for a process $p_i \in \calP_t$, the product
result is $u_{p_i} z_\text{local}$. After previous step, $z_\text{local}$ is
owned by $p_i$. Hence we only need to process the following
communication-free computation,
\begin{equation} \label{eq:prod-yuz}
    y_{p_i} = y_{p_i} + u_{p_i} z_\text{local},
\end{equation}
where $y_{p_i}$ is the output vector $y$ restricted to the subdomain in
$\Omega_t$ owned by $p_i$.

\bfemph{Dense submatrix.} Consider a dense submatrix associated with
$\Omega_t \times \Omega_s$ and process groups $\calP_t \times \calP_s$. If
there is only one process in $\calP_t$, then the matrix-vector
multiplication as in \eqref{eq:prod-dx} has already been conducted in the
first step and the result $z_\text{local}$ is also owned by $\calP_t$
after previous communication steps. Hence we simply add it to the output
vector,
\begin{equation} \label{eq:prod-yz}
    y_{\Omega_t} = y_{\Omega_t} + z_\text{local}.
\end{equation}
If there are more than one process in $\calP_t$, then the dense matrix is
owned by $\calP_t$ in a block row fashion and the matrix vector
multiplication admits,
\begin{equation}
    D_{\Omega_t \times \Omega_s} x_{\Omega_s} =
    \begin{pmatrix}
        d_{p_0} x_{\Omega_s} \\
        \vdots \\
        d_{p_{\lvert \calP_t \rvert -1}} x_{\Omega_s}
    \end{pmatrix}
    =
    \begin{pmatrix}
        d_{p_0} z_\text{local} \\
        \vdots \\
        d_{p_{\lvert \calP_t \rvert -1}} z_\text{local}
    \end{pmatrix},
\end{equation}
where $p_i \in \calP_t$ and each $p_i$ has a copy of $z_\text{local}$. In
this step, process $p_i$ is responsible for the following computation,
\begin{equation} \label{eq:prod-ydz}
    y_{p_i} = y_{p_i} + d_{p_i} z_\text{local},
\end{equation}
where $y_{p_i}$ is the same as that in \eqref{eq:prod-yuz}.

\begin{remark}
    Here we described the algorithm computing $y = \calK x$ for a
    distributed-memory \HMat{} $\calK$. A more standard matrix-vector
    multiplication operator in linear algebra would be $y = \alpha \calK x +
    \beta y$, which is the ``GEMV'' operation in level 2 BLAS. Such an
    operation can be easily adopted here if we do not initialize $y$ as a
    zero vector and modify \eqref{eq:prod-yuz}, \eqref{eq:prod-yz}, and
    \eqref{eq:prod-ydz} accordingly. All the rest steps remain unchanged.
\end{remark}

\subsection{Complexity Analysis}
\label{sec:Complexity}

In this section, we analyze the computational and the communication
complexities of the distributed-memory \HMat{}-vector multiplication
algorithm. To simplify the notation, we denote $L_p = O(\log P)$ and $L_N =
O( \log N )$ as the number of levels in process trees\footnote{Here we count
the number of levels in a process tree until the first level such that all
process groups contain one process.} and domain trees respectively.

The computational complexity is easy to conclude given our previous analysis
on the data balancing in Section~\ref{sec:data-distribution}. Notice that
our total number of floating-point operations stay identical to that of
sequential \HMat{}-vector multiplication if the extra computation in
tree-reduction is excluded. While, the computation in tree-reduction is of
lower order comparing to that of dense matrix-vector multiplication
conducted on each processes. Hence, the extra computation in communication
steps can be ignored in our complexity analysis. Further, processes conduct
float operations proportional to amounts of data they owned. Thanks to the
balanced data distribution, we conclude that the computational operations
are also balanced across all processes and each process conduct
$O\big(\frac{N \log N}{P} \big)$ operations.

The communication complexity consists of two parts: the latency ($\alpha$)
and the per-process inverse bandwidth ($\beta$). The complexity analysis
for the latency is relatively simpler and stay the same for different
admissibility conditions. The latency is essentially counting the number
of send/receive communications. Each process in the tree-reduction and
tree-broadcast steps conducts a reduction and broadcast among constant
number of processes. Hence each process conduct $O(1)$ send/receive
communications on each level. Summing all $L_P$ levels together, the
latencies for both tree-reduction and tree-broadcast are $O(\alpha \log
P)$. Regarding the Step 3 in our algorithm, as discussed in
Remark~\ref{rmk:step3}, each process only communicates with $O(\log P)$
other processes. Hence the latency for Step 4 and the overall latency are
$O(\alpha \log P)$. The complexities of inverse bandwidth, however, are
different for different admissibility conditions and are discussed
separately.

\bfemph{Weak admissibility condition.}
Consider the tree-reduction and tree-broadcast steps. At a given level $\ell
\leq L_P$, each process only participates one process group and owns a
constant number of submatrices. Hence the final concatenated array is of
length $O(L_P)$. Process 0 is the most communication intensive process. For
level $\ell = 1, \dots, L_P$, it communicates an array of size $O(\ell)$ in
both tree-reduction and tree-broadcast. Therefore, process 0 in total send
and receive $O(L_P^2)$ data, which is an upper bound for other processes.
The inverse bandwidth complexities for the tree-reduction and tree-broadcast
steps are then $O(\beta \log^2 P)$.

The inverse bandwidth complexity for the third step is very much simplified
for \HMats{} under weak admissibility condition due to one crucial
difference between weak admissibility condition and other admissibility
conditions. \HMats{} under weak admissibility condition only have
$\calH$-submatrices along their diagonal blocks, whereas \HMats{} under
other admissibility conditions have $\calH$-submatrices on off-diagonal
blocks. Under the distributed-memory setting, such a property means that the
source and target process groups remain the same for all $\calH$-submatrices
when weak admissibility condition is adopted. Hence only low-rank
submatrices are distributed among different source and target process
groups. Now we again consider process 0, who are group leaders across all
levels. For levels below $L_P$, process 0 does not participate any
submatrices with different source and target process groups. For level $L_P$
and above, process 0 is responsible to send the entire reduced array of
length $O(L_P)$ to other processes. Hence the inverse bandwidth complexities
for process 0 is $O(\beta \log P)$, which is the upper bound for other
processes.

Overall, the complexity, including both computational complexity and
communication complexity, for distributed-memory \HMats{} under weak
admissibility condition on $P$ processes is
\begin{equation} \label{eq:complexity-weak}
    O \bigg( \frac{N \log N}{P} + \alpha \log P + \beta \log^2 P \bigg).
\end{equation}

\bfemph{Standard admissibility condition.}
All communication complexity analyses under the weak admissibility condition
carry over to that under the standard admissibility condition with a
different prefactor, which is determined by the number of admissible
neighbors. Some extra communication costs come from those
$\calH$-submatrices singly owned by different target process and source
process. In this case, no tree-communication is needed. But the source
process need to pack all local data in this $\calH$-submatrices and send
them to the target process. The amount of local data in the
$\calH$-submatrices is a constant times the number of low-rank and dense
submatrices. Such $\calH$-submatrices are mostly corresponding to
neighboring subdomains and are of sizes $\frac{N}{P}$. With a complicated
calculation, which is omitted here, such $\calH$-submatrices have
$O\big(\log \frac{N}{P} \big)$ low-rank submatrices and
$O\Big(\big(\frac{N}{P} \big)^{\frac{d-1}{d}} \Big)$ dense submatrices,
where $d$ is the dimension of the problem. The number of low-rank
submatrices essentially calculates the number of levels whereas the number
of dense submatrices calculates the number of the subdomains of finest scale
on the interface of the two neighboring subdomains. Hence the extra
communication cost under standard admissibility condition is $O\Big( \beta
\Big( \log \frac{N}{P} + \big(\frac{N}{P} \big)^{\frac{d-1}{d}} \Big)
\Big)$.

Overall, the complexity for distributed-memory \HMats{} under standard
admissibility condition on $P$ processes is
\begin{equation} \label{eq:complexity-standard}
    O \bigg( \frac{N \log N}{P} + \alpha \log P + \beta \bigg( \log^2 P
    + \log \frac{N}{P} + \Big(\frac{N}{P} \Big)^{\frac{d-1}{d}} \bigg)
    \bigg).
\end{equation}

\begin{remark}
    According to \eqref{eq:complexity-weak} and
    \eqref{eq:complexity-standard}, we notice the trade-off between the
    computational complexity and the communication complexity. When $P$ is
    much smaller than $N$, the dominate cost comes from the computational
    part. While as $P$ approaches $N$, the computational cost is then
    $O(\log N)$ whereas the communication complexity is $O(\log^2 P)$
    dominating the cost.
\end{remark}

\section{Numerical Results}
\label{sec:numerical-results}

All numerical experiments were performed on the Texas Advanced Computing
Center (TACC) cluster, Stampede2. This cluster has $4,200$ Intel Knights
Landing nodes, each with 68 cores, 96~GB of DDR memory. Nodes are
interconnected via Intel Omni-Path network with a fat tree topology. We
allocate various number of nodes for our tests and each node runs 32 MPI
processes. The memory limit per process is 3~GB.

In the following numerical results, we adopt a few measurements to
demonstrate the parallel efficiency of our algorithm. In addition to the
regular wall-clock time (walltime), we also calculate the speedup as well as
the efficiency factor. Given a problem, we denote $P_0$ as the smallest
number of processes that are able to solve the problem and solve it in $t_0$
seconds. Meanwhile, solving the problem among $P_1$ processes for $P_1 \geq
P_0$ takes $t_1$ seconds. The speedup and the efficiency factor (percentage)
in this case are,
\begin{equation}
    \text{Speedup} = \frac{P_0 t_0}{t_1} \quad \text{ and} \quad
    \text{Eff} = \frac{P_0 t_0}{P_1 t_1} \cdot 100,
\end{equation}
respectively.

\subsection{\HMats{} for Two-Dimensional Problems}
\label{sec:numerical-results-2d}

Let $\Omega = [0,1]^2$ be the domain of interest. We discretize the problem
with $n$ points on each dimension for $n = 512, 1024, \dots, 65536$. Hence
the corresponding matrices are of size varying from $512^2 \times 512^2$ up
to $65536^2 \times 65536^2$. The structure of an \HMat{} is then determined
by a hierarchical partition of $\Omega$. Since the construction of \HMat{}
is beyond the scope of this paper and \HMat{}-vector multiplication does not
rely on the properties of the underlying problems, we fill dense submatrices
and low-rank submatrices in \HMats{} by random numbers and use these random
\HMats{} to explore the parallel scaling of our algorithm. Also random input
vectors are used in our tests. Both weak admissibility condition and
standard admissibility condition are explored. In addition, we use two
choices of $r$, $r=4$ and $r=8$, where the later makes problems more
computation intensive. Each \HMat{} is distributed among various number of
processes, from $32$ up to $16384$. The reported runtime is averaged over
$128$ random input vectors.

\begin{figure}[htp]
    \centering
    \includegraphics[width=0.4\textwidth]{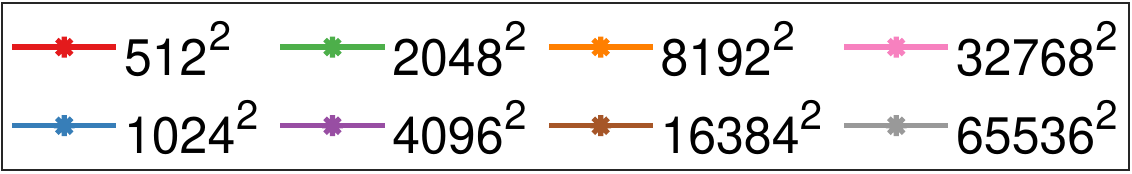}
    \vspace*{1em}

    \subfigure[Weak admissibility, $r=4$]{
        \includegraphics[width=0.45\textwidth]{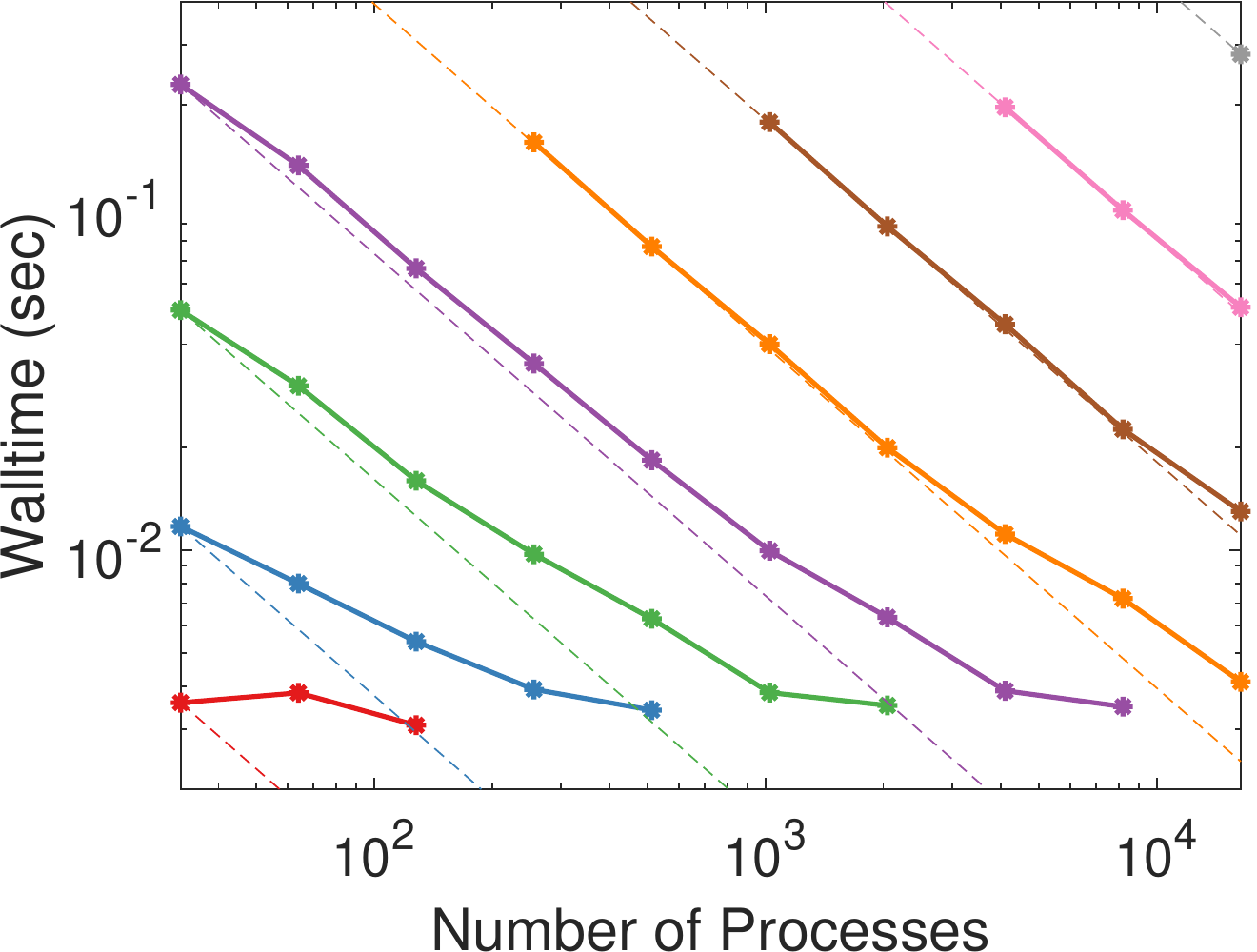}
        \label{fig:perf-weak-R4-2D-MV}}
    \quad
    \subfigure[Weak admissibility, $r=8$]{
        \includegraphics[width=0.45\textwidth]{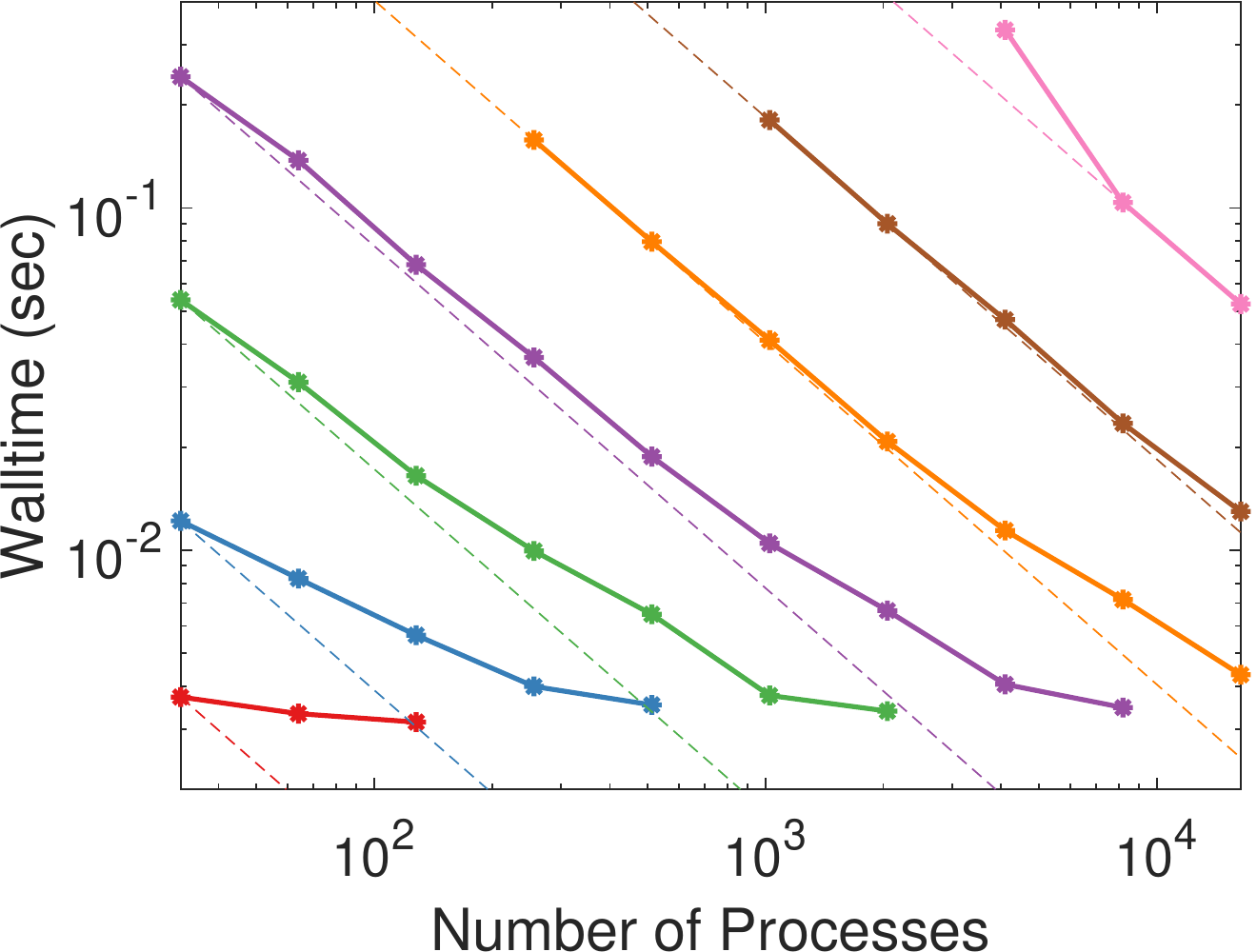}
        \label{fig:perf-weak-R8-2D-MV}}
    \vspace*{1em}

    \subfigure[Standard admissibility, $r=4$]{
        \includegraphics[width=0.45\textwidth]{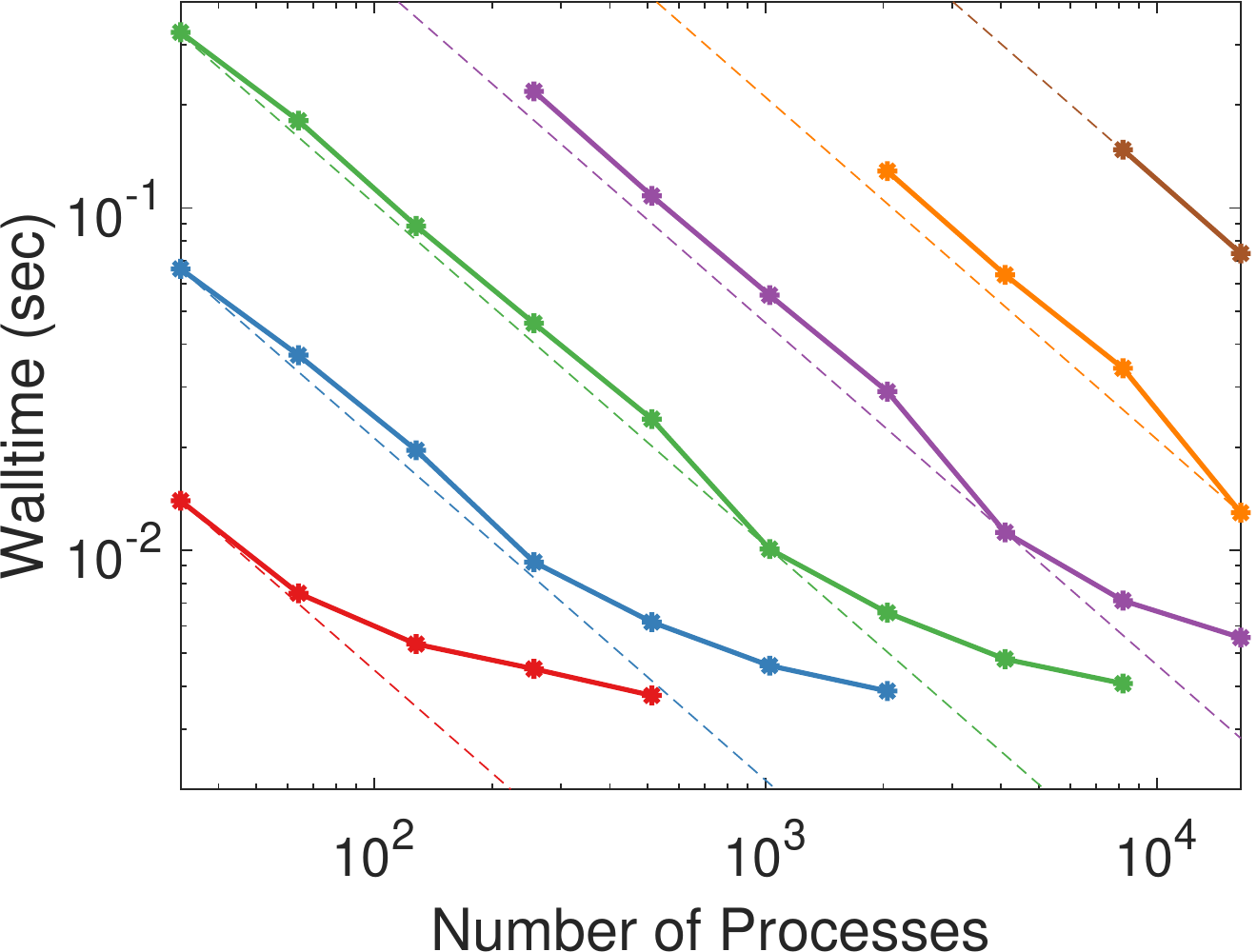}
        \label{fig:perf-standard-R4-2D-MV}}
    \quad
    \subfigure[Standard admissibility, $r=8$]{
        \includegraphics[width=0.45\textwidth]{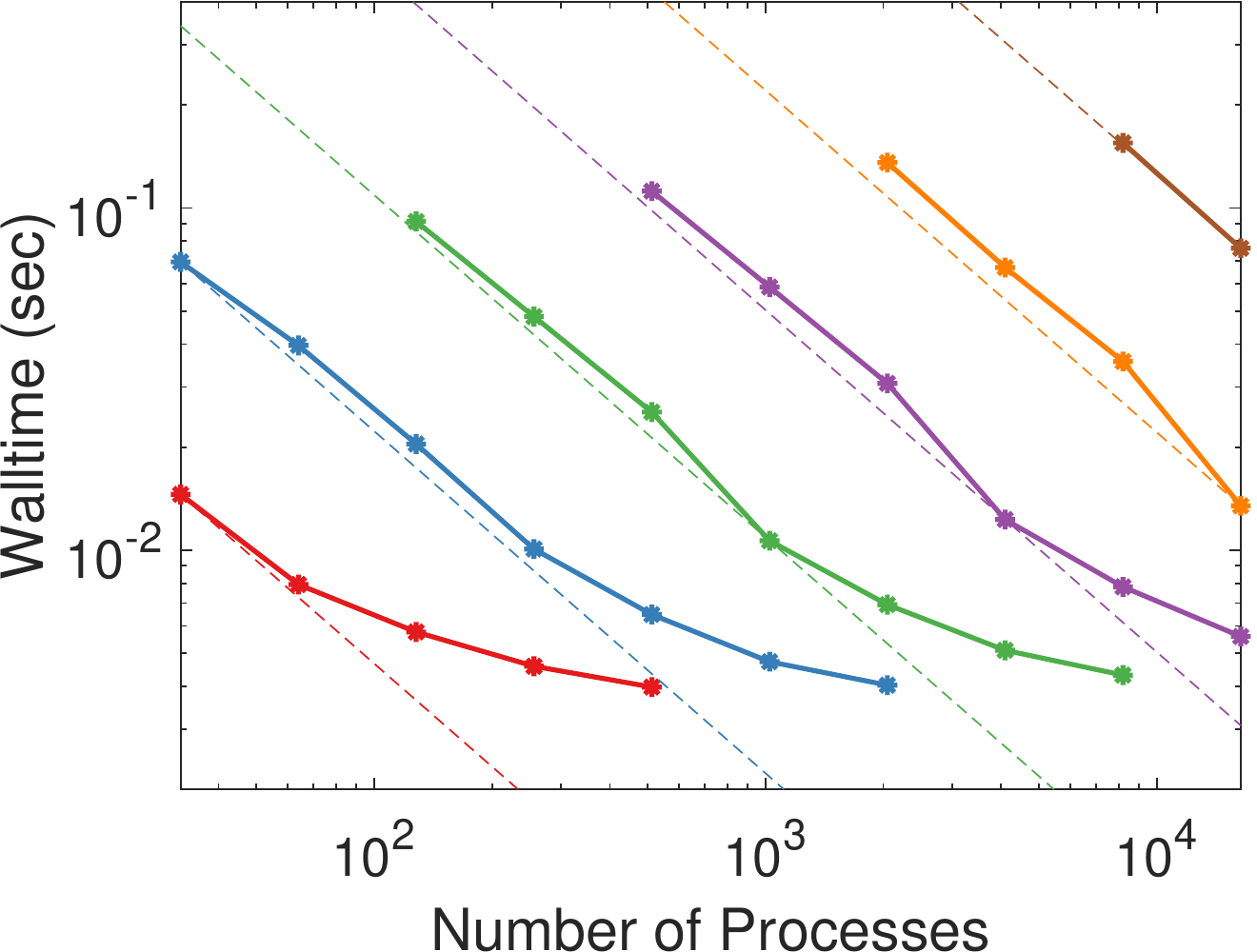}
        \label{fig:perf-standard-R8-2D-MV}}
        
    \caption{Strong scaling of \HMat{}-vector multiplication for various
    two-dimensional problems on various number of processes (up to 16384
    processes). Figure~(a) and (b) are \HMats{} under weak admissibility
    condition with rank being $4$ and $8$ respectively. Figure~(c) and (d)
    are \HMats{} under standard admissibility condition with rank being $4$
    and $8$ respectively. Solid lines are strong scaling curves and dash
    lines are their corresponding theoretical references. Different colors
    are problems of different sizes as indicated in the legend.}
    \label{fig:2d-results}
\end{figure}

\begin{table}
    \centering
    \begin{tabular}{llrrrrrrr}
\toprule
\multirow{2}{*}{$N$} & \multirow{2}{*}{$r$} & \multirow{2}{*}{$P$}
& \multicolumn{3}{c}{Weak} & \multicolumn{3}{c}{Standard} \\
\cmidrule(lr){4-6} \cmidrule(lr){7-9}
& & & Time (s) & Speedup & Eff (\%)
    & Time (s) & Speedup & Eff (\%) \\
\toprule
\multirow{3}{*}{$512^2$} & \multirow{3}{*}{4} &
      32 & 3.58e-03 &   32.0x & 100.0 & 1.39e-02 &   32.0x & 100.0 \\
& &   64 & 3.83e-03 &   30.0x &  46.8 & 7.48e-03 &   59.6x &  93.2 \\
& &  128 & 3.08e-03 &   37.2x &  29.1 & 5.31e-03 &   84.0x &  65.6 \\
\midrule
\multirow{5}{*}{$1024^2$} & \multirow{5}{*}{4} &
      32 & 1.17e-02 &   32.0x & 100.0 & 6.64e-02 &   32.0x & 100.0 \\
& &   64 & 7.98e-03 &   47.0x &  73.5 & 3.71e-02 &   57.2x &  89.4 \\
& &  128 & 5.40e-03 &   69.5x &  54.3 & 1.96e-02 &  108.5x &  84.8 \\
& &  256 & 3.91e-03 &   95.9x &  37.4 & 9.23e-03 &  230.3x &  90.0 \\
& &  512 & 3.40e-03 &  110.2x &  21.5 & 6.16e-03 &  345.0x &  67.4 \\
\midrule
\multirow{7}{*}{$2048^2$} & \multirow{7}{*}{4} &
      32 & 5.03e-02 &   32.0x & 100.0 & 3.26e-01 &   32.0x & 100.0 \\
& &   64 & 3.02e-02 &   53.3x &  83.2 & 1.80e-01 &   58.0x &  90.6 \\
& &  128 & 1.59e-02 &  100.9x &  78.8 & 8.86e-02 &  117.9x &  92.1 \\
& &  256 & 9.72e-03 &  165.5x &  64.6 & 4.60e-02 &  226.9x &  88.6 \\
& &  512 & 6.30e-03 &  255.4x &  49.9 & 2.42e-02 &  432.4x &  84.4 \\
& & 1024 & 3.83e-03 &  419.9x &  41.0 & 1.01e-02 & 1037.1x & 101.3 \\
& & 2048 & 3.51e-03 &  457.8x &  22.4 & 6.56e-03 & 1592.4x &  77.8 \\
\midrule
\multirow{9}{*}{$4096^2$} & \multirow{9}{*}{4} &
      32 & 2.30e-01 &   32.0x & 100.0 &        - &       - &     - \\
& &   64 & 1.33e-01 &   55.2x &  86.2 &        - &       - &     - \\
& &  128 & 6.66e-02 &  110.4x &  86.2 &        - &       - &     - \\
& &  256 & 3.51e-02 &  209.4x &  81.8 & 2.19e-01 &  256.0x & 100.0 \\
& &  512 & 1.83e-02 &  401.3x &  78.4 & 1.08e-01 &  517.7x & 101.1 \\
& & 1024 & 9.95e-03 &  738.6x &  72.1 & 5.56e-02 & 1009.2x &  98.6 \\
& & 2048 & 6.36e-03 & 1155.7x &  56.4 & 2.91e-02 & 1932.2x &  94.3 \\
\midrule
\multirow{7}{*}{$8192^2$} & \multirow{7}{*}{4} &
     256 & 1.56e-01 &  256.0x & 100.0 &        - &       - &     - \\
& &  512 & 7.71e-02 &  516.5x & 100.9 &        - &       - &     - \\
& & 1024 & 4.00e-02 &  995.6x &  97.2 &        - &       - &     - \\
& & 2048 & 1.99e-02 & 1999.3x &  97.6 & 1.28e-01 & 2048.0x & 100.0 \\
& & 4096 & 1.11e-02 & 3577.2x &  87.3 & 6.38e-02 & 4118.2x & 100.5 \\
& & 8192 & 7.21e-03 & 5523.2x &  67.4 & 3.40e-02 & 7723.1x &  94.3 \\
& & 16384 & 4.12e-03 & 9666.5x &  59.0 & 1.29e-02 & 20411.5x & 124.6 \\
\midrule
\multirow{5}{*}{$16384^2$} & \multirow{5}{*}{4} &
    1024 & 1.78e-01 & 1024.0x & 100.0 &        - &       - &     - \\
& & 2048 & 8.83e-02 & 2066.0x & 100.9 &        - &       - &     - \\
& & 4096 & 4.57e-02 & 3991.3x &  97.4 &        - &       - &     - \\
& & 8192 & 2.26e-02 & 8091.8x &  98.8 & 1.48e-01 & 8192.0x & 100.0 \\
& & 16384 & 1.30e-02 & 14063.9x &  85.8 & 7.36e-02 & 16471.7x & 100.5 \\
\midrule
\multirow{3}{*}{$32768^2$} & \multirow{3}{*}{4} &
    4096 & 1.97e-01 & 4096.0x & 100.0 &        - &       - &     - \\
& & 8192 & 9.86e-02 & 8188.2x & 100.0 &        - &       - &     - \\
& & 16384 & 5.13e-02 & 15744.5x &  96.1 &        - &       - &     - \\
\midrule
\multirow{1}{*}{$65536^2$} & \multirow{1}{*}{4} &
    16384 & 2.82e-01 & 16384.0x & 100.0 &        - &       - &     - \\
\bottomrule
    \end{tabular}
    \caption{Numerical results of distributed-memory \HMat{}-vector
    multiplication for two-dimensional problems.}
    \label{tab:num-result-2d}
\end{table}

Figure~\ref{fig:2d-results} depicts strong scaling plots for different
\HMats{} and Table~\ref{tab:num-result-2d} further details walltimes,
speedups and efficiency factors. In both weak admissibility condition cases,
Figure~\ref{fig:perf-weak-R4-2D-MV} and Figure~\ref{fig:perf-weak-R8-2D-MV},
strong scaling is well-preserved as we keep doubling the number of
processes. Towards the end of each curve, when the communication cost
dominates the walltime, the walltime remain flat for a long time, which
means that the communication cost grows very mildly as the number of
processes increases. In standard admissibility condition cases,
Figure~\ref{fig:perf-standard-R4-2D-MV} and
Figure~\ref{fig:perf-standard-R8-2D-MV}, good strong scaling is also
observed in most cases. Comparing to the weak admissibility condition cases,
especially towards the end of each curve, the communication cost kicks in
earlier as the number of processes increase, which is due to the different
prefactors in the complexity analysis in Section~\ref{sec:Complexity}.
Table~\ref{tab:num-result-2d} provides more evidences supporting our
comments. We emphasize that the parallel efficiencies are impressive
especially for larger problems. For example, in both $N=4096^2$ and
$N=8192^2$ cases, parallel efficiencies are above 72 percent in weak
admissibility condition cases and above 90 percent in standard admissibility
condition cases, even when thousands of processes are used. Finally, we
would like to comment on the weak scaling. Although not been plotted in
figures, weak scaling\footnote{The computational cost grows quasi-linearly
whereas the number of processes grows linearly. Here our weak scaling
definition ignores the extra logarithmic factor.} can be read from
connecting dots vertically in figures. Clearly, on the top half of each
figure, the weak scaling is near ideal (flat). Hence we claim that our
algorithm and implementation give numerical results of both good strong
scaling and weak scaling.

\subsection{\HMats{} for Three-Dimensional Problems}
\label{sec:numerical-results-3d}

In this section, we perform numerical results for domain $\Omega = [0,1]^3$.
We discretize the problem with $n$ being $64, 128, \dots, 1024$ and the
corresponding matrices are of size varying from $64^3 \times 64^3$ up to
$1024^3 \times 1024^3$. Similar as in the two-dimensional cases, we adopt
random \HMats{} and random input vectors to explore the parallel scaling of
our algorithm. Both weak admissibility condition and standard admissibility
condition are explored as well as two choices of $r$. Each \HMat{} is
distributed among various number of processes, from $32$ up to $16384$.
Reported runtime is averaged over $128$ random input vectors.

\begin{figure}[htp]
    \centering
    \includegraphics[width=0.3\textwidth]{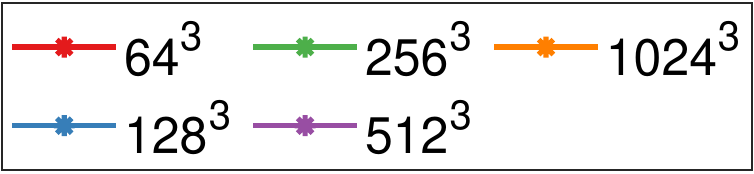}
    \vspace*{1em}

    \subfigure[Weak admissibility, $r=4$]{
        \includegraphics[width=0.45\textwidth]{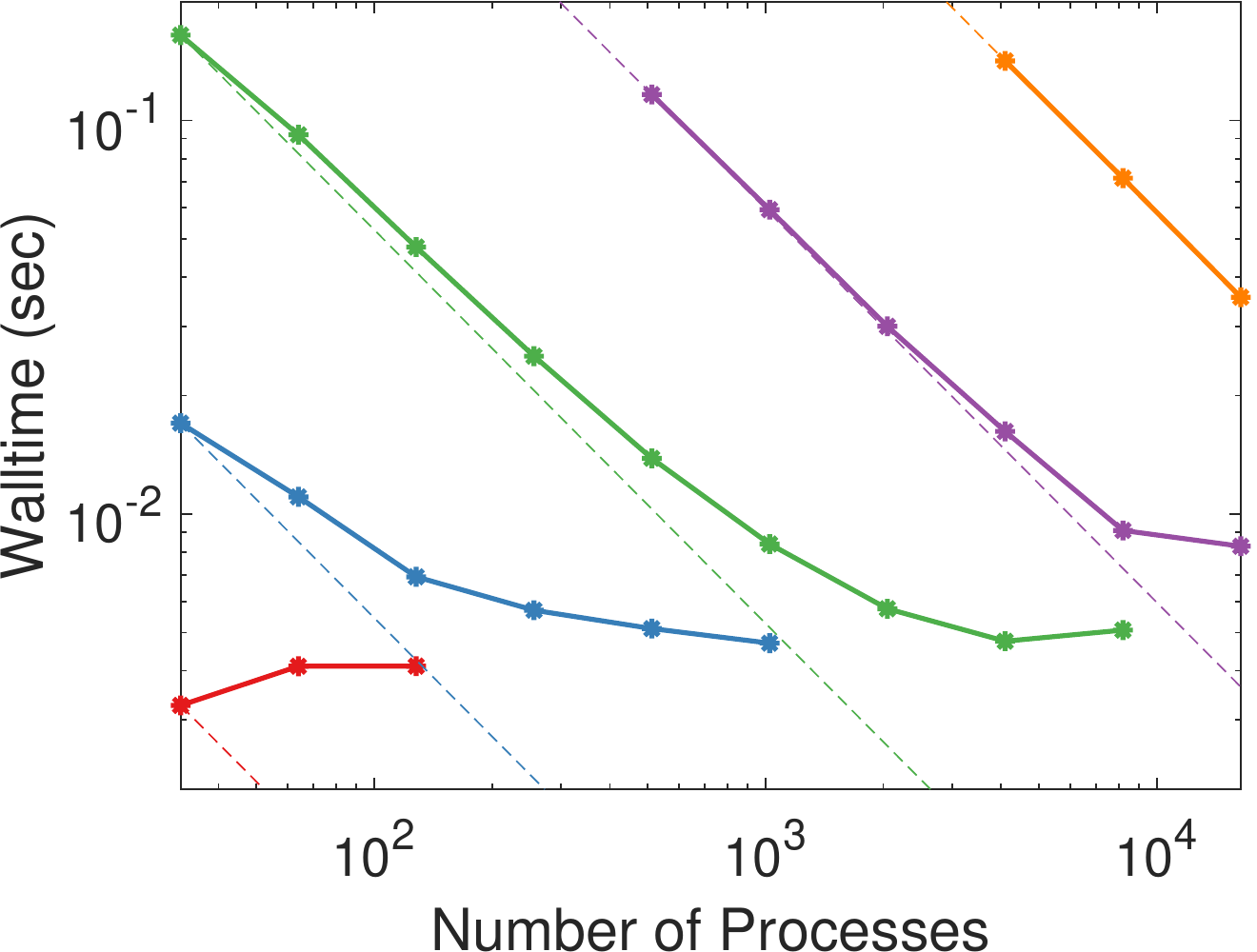}
        \label{fig:perf-weak-R4-3D-MV}}
    \quad
    \subfigure[Weak admissibility, $r=8$]{
        \includegraphics[width=0.45\textwidth]{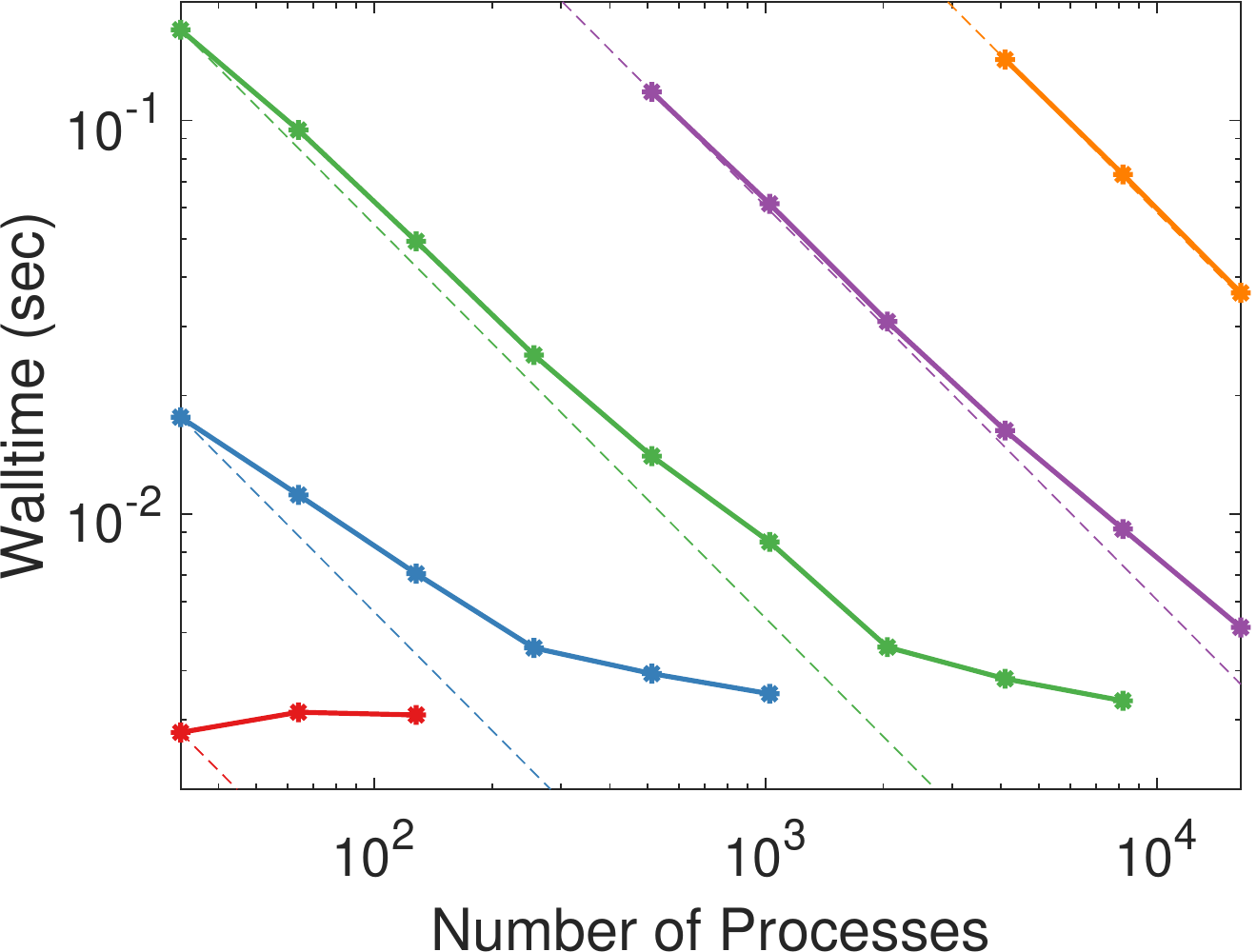}
        \label{fig:perf-weak-R8-3D-MV}}
    \vspace*{1em}

    \subfigure[Standard admissibility, $r=4$]{
        \includegraphics[width=0.45\textwidth]{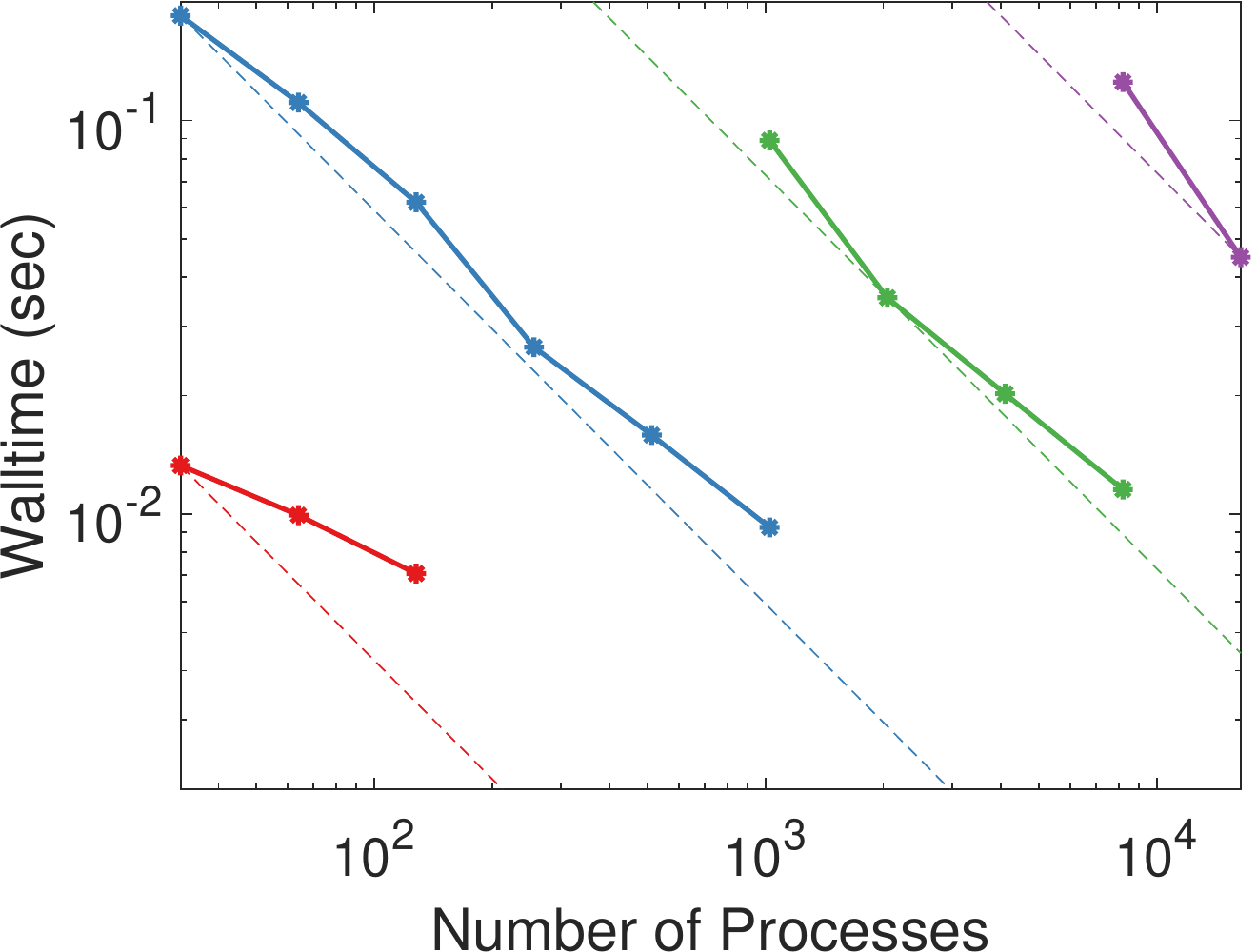}
        \label{fig:perf-standard-R4-3D-MV}}
    \quad
    \subfigure[Standard admissibility, $r=8$]{
        \includegraphics[width=0.45\textwidth]{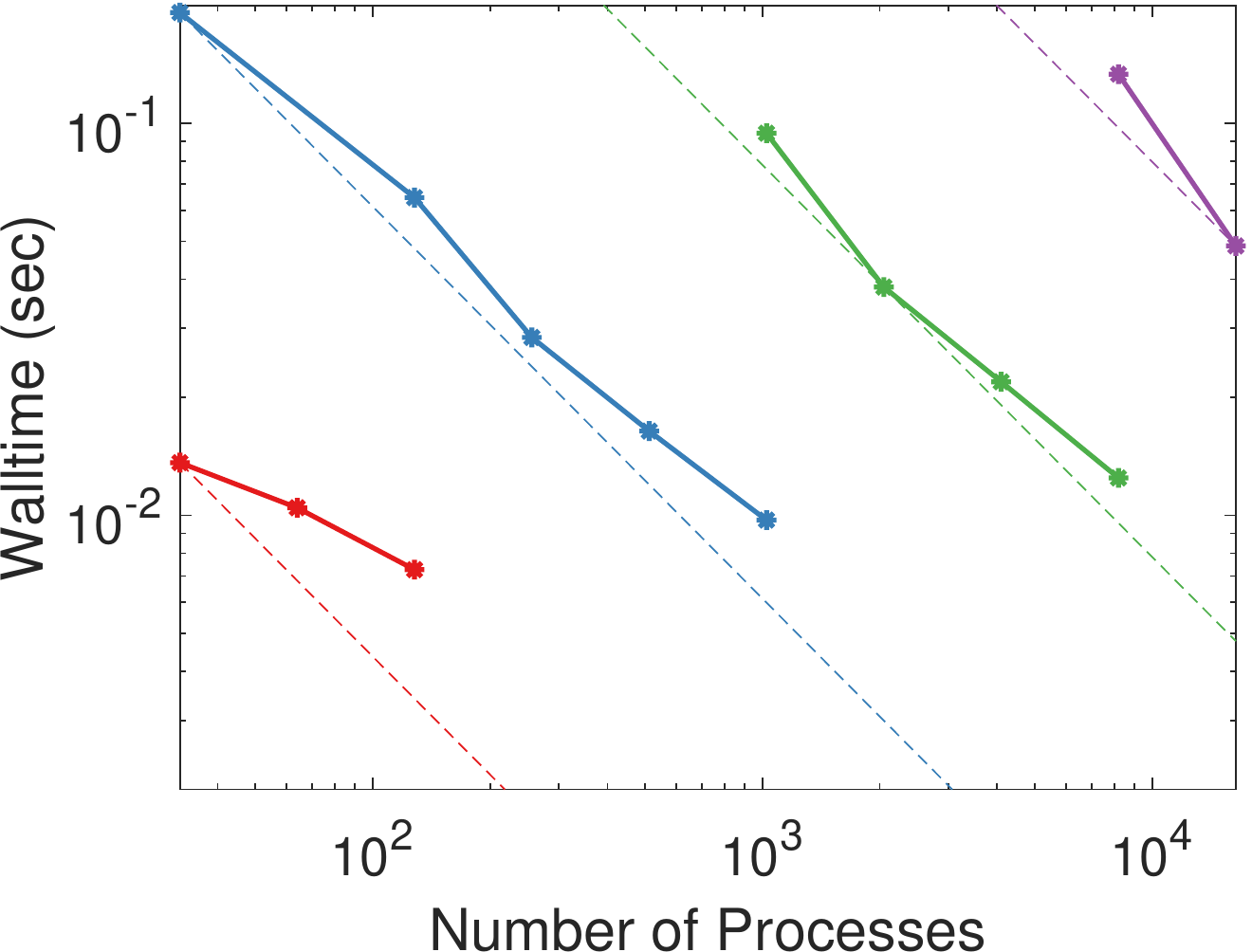}
        \label{fig:perf-standard-R8-3D-MV}}
        
    \caption{Strong scaling of \HMat{}-vector multiplication for various
    three-dimensional problems on various number of processes (up to 16384
    processes). Figure~(a) and (b) are \HMats{} under weak admissibility
    condition with rank being $4$ and $8$ respectively. Figure~(c) and (d)
    are \HMats{} under standard admissibility condition with rank being $4$
    and $8$ respectively. Solid lines are strong scaling curves and dash
    lines are their corresponding theoretical references. Different colors
    are problems of different sizes as indicated in the legend.}
    \label{fig:2d-results}
\end{figure}

\begin{table}
    \centering
    \begin{tabular}{llrrrrrrr}
\toprule
\multirow{2}{*}{$N$} & \multirow{2}{*}{$r$} & \multirow{2}{*}{$P$}
& \multicolumn{3}{c}{Weak} & \multicolumn{3}{c}{Standard} \\
\cmidrule(lr){4-6} \cmidrule(lr){7-9}
& & & Time (s) & Speedup & Eff (\%)
    & Time (s) & Speedup & Eff (\%) \\
\toprule
\multirow{3}{*}{$64^3$} & \multirow{3}{*}{4} &
      32 & 3.27e-03 &   32.0x & 100.0 & 1.33e-02 &   32.0x & 100.0 \\
& &   64 & 4.11e-03 &   25.5x &  39.8 & 9.94e-03 &   42.8x &  66.8 \\
& &  128 & 4.11e-03 &   25.4x &  19.9 & 7.06e-03 &   60.2x &  47.0 \\
\midrule
\multirow{6}{*}{$128^3$} & \multirow{6}{*}{4} &
      32 & 1.70e-02 &   32.0x & 100.0 & 1.85e-01 &   32.0x & 100.0 \\
& &   64 & 1.11e-02 &   49.2x &  76.9 & 1.11e-01 &   53.2x &  83.2 \\
& &  128 & 6.93e-03 &   78.6x &  61.4 & 6.20e-02 &   95.4x &  74.6 \\
& &  256 & 5.70e-03 &   95.6x &  37.3 & 2.66e-02 &  222.7x &  87.0 \\
& &  512 & 5.12e-03 &  106.4x &  20.8 & 1.59e-02 &  372.5x &  72.8 \\
& & 1024 & 4.70e-03 &  115.9x &  11.3 & 9.26e-03 &  638.9x &  62.4 \\
\midrule
\multirow{9}{*}{$256^3$} & \multirow{9}{*}{4} &
      32 & 1.65e-01 &   32.0x & 100.0 &        - &       - &     - \\
& &   64 & 9.20e-02 &   57.3x &  89.6 &        - &       - &     - \\
& &  128 & 4.77e-02 &  110.6x &  86.4 &        - &       - &     - \\
& &  256 & 2.52e-02 &  209.4x &  81.8 &        - &       - &     - \\
& &  512 & 1.39e-02 &  380.7x &  74.4 &        - &       - &     - \\
& & 1024 & 8.40e-03 &  628.4x &  61.4 & 8.91e-02 & 1024.0x & 100.0 \\
& & 2048 & 5.75e-03 &  917.6x &  44.8 & 3.55e-02 & 2572.3x & 125.6 \\
& & 4096 & 4.75e-03 & 1109.8x &  27.1 & 2.02e-02 & 4510.0x & 110.1 \\
& & 8192 & 5.07e-03 & 1039.8x &  12.7 & 1.16e-02 & 7894.9x &  96.4 \\
\midrule
\multirow{6}{*}{$512^3$} & \multirow{6}{*}{4} &
     512 & 1.16e-01 &  512.0x & 100.0 &        - &       - &     - \\
& & 1024 & 5.94e-02 & 1004.2x &  98.1 &        - &       - &     - \\
& & 2048 & 3.00e-02 & 1987.9x &  97.1 &        - &       - &     - \\
& & 4096 & 1.62e-02 & 3674.9x &  89.7 &        - &       - &     - \\
& & 8192 & 9.09e-03 & 6561.1x &  80.1 & 1.25e-01 & 8192.0x & 100.0 \\
& & 16384 & 8.29e-03 & 7196.6x &  43.9 & 4.50e-02 & 22811.5x & 139.2 \\
\midrule
\multirow{3}{*}{$1024^3$} & \multirow{3}{*}{4} &
    4096 & 1.42e-01 & 4096.0x & 100.0 &        - &       - &     - \\
& & 8192 & 7.14e-02 & 8131.5x &  99.3 &        - &       - &     - \\
& & 16384 & 3.56e-02 & 16320.9x &  99.6 &        - &       - &     - \\
\bottomrule
    \end{tabular}
    \caption{Numerical results of distributed-memory \HMat{}-vector
    multiplication for three-dimensional problems.}
    \label{tab:num-result-3d}
\end{table}

Comments for two-dimensional problemss as in
Section~\ref{sec:numerical-results-2d} apply seamless to three-dimensional
problems. Both under weak and standard admissibility condition cases, strong
scaling and weak scaling are well-preserved as the number of processes
increases. \HMats{} under weak admissibility condition show better parallel
efficiencies comparing to that under standard admissibility condition. Now
we focus on the comparison of two-dimensional problems and three-dimensional
problems. Comparing Figure~\ref{fig:perf-weak-R4-2D-MV} and
Figure~\ref{fig:perf-weak-R8-2D-MV} to Figure~\ref{fig:perf-weak-R4-3D-MV}
and Figure~\ref{fig:perf-weak-R8-3D-MV} respectively, we find that all four
figures show similar strong scaling as well as weak scaling. This behavior
has already been predicted by \eqref{eq:complexity-weak}, where the
complexity under weak admissibility condition is independent of the
dimensionality of the problem. While, comparing
Figure~\ref{fig:perf-standard-R4-2D-MV} and
Figure~\ref{fig:perf-standard-R8-2D-MV} to
Figure~\ref{fig:perf-standard-R4-3D-MV} and
Figure~\ref{fig:perf-standard-R8-3D-MV} respectively, two-dimensional
problems show better strong scaling than their three-dimensional
counterparts. Under standard admissibility condition, the number of
neighboring subdomains increases as the dimension increases, which also
implies that the required communication cost will increase. As detailed in
\eqref{eq:complexity-standard}, the communication complexity depends
monotonically on the dimension $d$. Hence, as proved by numerical resutls,
the communication cost dominate the walltime earlier for biger $d$.

\section{Conclusion}
\label{sec:conclusion}

In this paper, we introduce the data distribution of distributed \HMats{}
and a distributed-memory \HMats{}-vector multiplication algorithm.

Given the tree structure of the domain organization in \HMat{}, we also
organize our processes in a process tree. Two process trees are adopted for
the target and source domains. Under our data distribution scheme, the load
balancing factors are constants for both weak admissibility condition (the
constant is independent of dimension $d$) and standard admissibility
condition (the constant depends on $d$). For problems of extremely large
size $N$, our data distribution scheme allows the number of processes to
grow as big as $O(N)$. In this case, each process owns a part of the
\HMat{}, whose size depends only logarithmically on $N$. Therefore, our data
distribution is feasible for problems of extremely large sizes on massive
number of processes.

The proposed distributed-memory \HMat{}-vector multiplication algorithm is
parallel efficient. Specifically under our tree organizations of both
processes and data, we introduce a tree communication scheme, i.e.,
``tree-reduce'' and ``tree-broadcast'', to significantly reduce the latency
complexity. All required computations in sequential \HMat{}-vector
multiplication are evenly distributed among all processes. Importantly, our
algorithm totally avoids the expensive scheduling step, which is as
expensive as $\Omega(P^2)$ on $P$ processes. Overall, our algorithm
complexities for a $d$-dimensional problem of size $N$ distributed among $P$
processes are $O \Big( \frac{N \log N}{P} + \alpha \log P + \beta \log^2 P
\Big)$ and $O \Big( \frac{N \log N}{P} + \alpha \log P + \beta \Big( \log^2
P + \log \frac{N}{P} + \big(\frac{N}{P} \big)^{\frac{d-1}{d}} \Big) \Big)$
for weakly admissibility condition and standard admissibility condition
respectively, where $\alpha$ denotes the latency and $\beta$ denotes the
per-process inverse bandwidth.

There are several future directions for improvement, both in algorithm and
in implementation. Instead of pure ``MPI'' parallelization, one can combine
``OpenMP'' and ``MPI'' to further reduce the local communications within a
node. This could improve the communication complexity, especially for
\HMats{} under standard admissibility condition, by a big factor. Other
\HMat{} algebraic operations can also be efficiently parallelized given our
data distribution and process organization. In a companion paper, we will
introduce distributed-memory \HMat{} compression, \HMat{} addition, as well
as \HMat{}-\HMat{} multiplication.

{\bf Availability.} The distributed-memory \HMat{} code, DMHM, is available
under the GPLv3 license at \url{https://github.com/YingzhouLi/dmhm}. The
code support both two-dimensional and three-dimensional problems.

\section*{Acknowledgments}

The authors acknowledge the Texas Advanced Computing Center (TACC) at The
University of Texas at Austin for providing HPC resources that have
contributed to the research results reported within this paper. The work of
Y.L. is supported in part by the US National Science Foundation under awards
DMS-1454939 and DMS-2012286, and by the US Department of Energy via grant
DE-SC0019449. The work of L.Y. is partially supported by the U.S. Department
of Energy, Office of Science, Office of Advanced Scientific Computing
Research, Scientific Discovery through Advanced Computing (SciDAC) program
and the National Science Foundation under award DMS-1818449.

\bibliographystyle{chicago}
\bibliography{dmhm}

\end{document}